\documentclass[a4paper]{amsart}
\usepackage[latin9]{inputenc}
\usepackage{a4}
\usepackage{amsfonts}
\usepackage{amssymb}
\usepackage{amsmath}
\usepackage{amsthm}
\usepackage{changepage}
\usepackage{nomencl}
\usepackage{geometry}
\usepackage{tikz}
\usepackage{verbatim}
\usetikzlibrary{matrix}
\usepackage{hyperref}
\usepackage[toc,page]{appendix}
\usepackage{mathrsfs}
\usepackage[square, numbers, comma]{natbib}

\usepackage{enumerate}
\usepackage{amsbsy}

\usepackage[all]{xy}
\usepackage{pb-diagram,pb-xy}
\bibpunct{[}{]}{,}{n}{}{,}
\input cyracc.def

\begin{document}
	\providecommand{\keywords}[1]{\textbf{\textit{Keywords: }} #1}
	\newtheorem{theorem}{Theorem}[section]
	\newtheorem{lemma}[theorem]{Lemma}
	\newtheorem{proposition}[theorem]{Proposition}
	\newtheorem{corollary}[theorem]{Corollary}
	\newtheorem{problem}[theorem]{Problem}
	\newtheorem{question}[theorem]{Question}
	\newtheorem{conjecture}[theorem]{Conjecture}
	\newtheorem{claim}[theorem]{Claim}
	\newtheorem{condition}[theorem]{Condition}

	\theoremstyle{definition}
	\newtheorem{definition}[theorem]{Definition} 
	\theoremstyle{remark}
	\newtheorem{remark}[theorem]{Remark}
	\newtheorem{example}[theorem]{Example}
	\newtheorem{condenum}{Condition}
	
	\def\p{\mathfrak{p}}
	\def\q{\mathfrak{q}}
	\def\s{\mathfrak{S}}
	\def\Gal{\mathrm{Gal}}
	\def\Ker{\mathrm{Ker}}
	\def\soc{\mathrm{soc}}
	\def\Coker{\mathrm{Coker}}
	\newcommand{\cc}{{\mathbb{C}}}   
	\newcommand{\ff}{{\mathbb{F}}}  
	\newcommand{\nn}{{\mathbb{N}}}   
	\newcommand{\qq}{{\mathbb{Q}}}  
	\newcommand{\rr}{{\mathbb{R}}}   
	\newcommand{\zz}{{\mathbb{Z}}}  
	\def\K{\kappa}
	
	\title[Reducible fibers of polynomials over integers]{Reducible fibers of polynomials over integers: the Hilbert-Siegel problem over $\mathbb{Q}$ and its applications}
	\author{Joachim K\"onig}
	\address{Department of Mathematics Education, Korea National University of Education, Cheongju 28173, South Korea}
	\email{jkoenig@knue.ac.kr}
	\begin{abstract}
		Given a polynomial $f\in \mathbb{Q}[X]$, the Hilbert--Siegel problem asks for a description of the set of all values $t_0\in \mathbb{Z}$, possibly up to a finite set, such that $f(X)-t_0$ is reducible. Reducibility occurs trivially for $t_0$ in the value set of some non-linear left-composition factor of $f$. The problem therefore reduces to determining the polynomials having infinitely many further reducible values. The biggest step towards a solution so far was taken in \cite{BKN26}. We strengthen some of the methods obtained there to present the complete solution of the Hilbert--Siegel problem over $\mathbb{Q}$.
	We give applications to the problem of newly reducible iterates of polynomials, and to an unlikely intersection problem on simultaneous reducible fibers of two polynomial maps over $\mathbb{Q}$, generalizing a result by Bilu and Tichy.
	\end{abstract}
	\maketitle

	\section{Introduction and main results}
	Let $f\in \mathbb{Z}[X]$ be a nonconstant integer polynomial. Hilbert's irreducibility theorem (applied to the irreducible two-variable polynomial $f(X)-t$) shows that $f(X)-t_0$ is irreducible over $\mathbb{Q}$ for all integer values $t_0$ outside of some set of relative density $0$. On the other hand, the set $\mathcal{R}_f:=\{t_0\in \mathbb{Z}\mid f(X)-t_0 \textrm{ reducible over } \mathbb{Q}\}$ is of course infinite as soon as $\deg(f)>1$: Indeed, any value $t_0 \in f(\mathbb{Q})\cap \mathbb{Z}$, and more generally every value $t_0\in f_1(\mathbb{Q})\cap \mathbb{Z}$, where $f_1\in \mathbb{Q}[X]$ is a non-linear left-compositional factor of $f$ (i.e., $f=f_1\circ f_2$ for some $f_2\in \mathbb{Q}[X]$) belongs to the set $\mathcal{R}_f$. These may be viewed as ``trivially reducible values" (although the explicit determination of all left-compositional factors of $f$ is itself a nontrivial computational problem, related to Ritt's theorems on polynomial decompositions).  The {\it Hilbert-Siegel problem}\footnote{The name is due to the combination of Hilbert's irreducibility theorem and Siegel's theorem on finiteness of integral points needed for the solution. Of course the problem is equally relevant over rings of integers of more general number fields (and the original formulation of the problem does account for this).}, proposed in \cite{Fried} (see also \cite{Debes99}), asks for the full description of the set $\mathcal{R}_f$ up to finitely many values, and in particular for the classification of polynomials $f$ for which $\mathcal{R}_f$ contains infinitely many elements other than the trivial ones described above. 
   The polynomial $f=X^n$ is such an example as soon as $4|n$. Indeed, by a classical result known as the Vahlen-Capelli theorem \cite{Vahlen}, a necessary and sufficient condition for $X^n-a\in \mathbb{Z}[X]$ to be reducible is that $a$ is a $p$-th power for some prime divisor $p$ of $n$ {\it or} ($4|n$ and $-4a$ is a fourth power). The former alternative accounts for the compositional left-factors $X^p$ of $X^n$; however $-\frac{1}{4}X^4$ does not have a common left-factor with $X^n$ since $-1$ is not a square.
	
	Via Siegel's theorem \cite{Siegel}, the Hilbert-Siegel problem reduces to the determination of {\it Siegel functions} $g\in \mathbb{Q}(X)$ such that (the numerator of) $f(X)-g(Y)$ is reducible. Monodromy allows for a group-theoretical description of polynomials as well as Siegel functions in terms of branch cycle descriptions, turning the Hilbert-Siegel problem into a permutation-group theoretical problem, albeit a very intricate one. Fried \cite{Fried74}, \cite{Fried} first accepted the group-theoretical challenges, leading to a solution in the case of {\it indecomposable} polynomials 
	  $f$ (corresponding to {\it primitive} monodromy groups), but the general case is still considered much harder. The biggest leap yet towards a full solution was undertaken in \cite{BKN26}, solving the problem in full for polynomials not factoring through an indecomposable of degree $2$ or $4$, and in fact reducing the remaining problem (over $\mathbb{Q}$) to the determination of {\it minimally reducible pairs} $(f,g)$ only in the special (although admittedly difficult) case where $\textrm{Mon}_{\mathbb{Q}}(f) = \textrm{Mon}_{\mathbb{Q}}(g)$ is solvable, $f$ is of even degree, $g$ is of degree $2\deg(f)$, and $g$ is a ``non-polynomial" Siegel function over $\mathbb{Q}$, i.e., a Siegel function with $|g^{-1}(\infty)|=2$ and both poles having the same order.

	In this article, we give a full description of all such pairs. Concretely, our main contribution is:
	\begin{theorem}
		\label{thm:main}
		Let $(f,g)$ be a minimally reducible pair, where $f\in \mathbb{Q}[X]$ is a polynomial with solvable monodromy group, and $g\in \mathbb{Q}(X)$ is a Siegel function over $\mathbb{Q}$ with $|g^{-1}(\infty)|=2$. Then $\deg(f)=8$, $\deg(g)=16$.\\
		Conversely, there are three explicitly computable families of minimally reducible pairs $(f,g)$ of those degrees.
		\end{theorem}
		Parameterizations of the three families will be given in Section \ref{sec:ex}.
	
	Together with the classification of minimally reducible pairs $(f,g)$ such that either $g$ is also a polynomial or $f$ has nonsolvable monodromy group, obtained in \cite{BKN26}, and with some explicit computations, this leads to a full solution of the Hilbert--Siegel problem over $\mathbb{Q}$ in the following form. To avoid misunderstanding, it should be emphasized that the below result along with all of its applications does not replace, but rather crucially requires the results of \cite{BKN26}.
	\begin{corollary}
		\label{cor:hilbert_siegel}
		The set $\mathcal{R}_f$ of reducible integer values of a polynomial $f\in \mathbb{Q}[X]$ is contained in the union of a finite set, of the value sets $h(\mathbb{Q})\cap \mathbb{Z}$ of indecomposable left-composition factors $h$ of $f$, and the value sets $\varphi_i(\mathbb{Q})\cap \mathbb{Z}$ of $j\ge 0$ further Siegel functions $\varphi_1,\dots,\varphi_j \in\mathbb{Q}(X)$ such that one of the following holds:
		\begin{itemize}
			\item[a)] $j=0$, i.e., there are no extra Siegel functions needed to describe $\mathcal{R}_f$.
			\item[b)] $j=1$; and either $f$ has a left-composition factor $h:=f_1\circ f_2$ where $f_1,f_2\in \mathbb{Q}[X]$ are indecomposable of degrees $(\deg(f_1),\deg(f_2))\in \{(2,2), (4,2), (2,4)\}$, where the cases involving a quartic are as in Theorems \ref{thm:8t44} and \ref{thm:8t47} respectively, up to composition with linears; or $f$ has a nonsolvable left-composition factor  $f_1 = \mu\circ  X(X^2+5aX+b)^2$ of degree $5$,\footnote{More precisely, $a,b\in \mathbb{Q}$ need to be such that  $b=5a^2+y^2-5x^2$ for suitable $(x,y)\in \mathbb{Q}^2\setminus\{(0,0)\}$. This family is known (in a different, but equivalent parameterization, only seemingly more complicated due to taking into account the composition with linears) from \cite[Section 2.3]{DF}.} $\mu\in \mathbb{Q}[X]$ linear.
			\item[c)] $j=2$; and either $f$ has a left-composition factor $\mu\circ a^2X^4(aX^2+4)^2$, for $\mu\in \mathbb{Q}[X]$ linear and $a\in \mathbb{Q}^\times$ \footnote{More precisely, $a$ needs to be of the form $a=2u^2-v^2$, see Theorem \ref{thm:8t28}.}, or  $f$ has a left-composition factor $\mu\circ X^4(X^4+c)^4$ for $\mu\in \mathbb{Q}[X]$ linear and $c\in \mathbb{Q}^\times$.
			\end{itemize}
		\end{corollary}
		
		Previously it was not even clear whether there is any absolute bound $j\in \mathbb{N}$ (independent of the polynomial $f$) such that the nontrivial reducible values can always be parameterized by at most $j$ functions.
		Note that in all cases, the Siegel functions $\varphi_i$ are explicitly computable and depend only on the exhibited left-composition factors. To fully answer the original finiteness question of the Hilbert-Siegel problem, it remains to figure out, for the cases in b), c), when precisely the associated Siegel functions $\varphi_i$ have infinitely many integral values. Since, under the given specifications of parameters, they automatically have infinitely many rational values of bounded denominator, infiniteness is certainly reached upon manipulating the linear left-factor $\mu$ accordingly, and in particular holds for suitable {\it integer} polynomials $f\in \mathbb{Z}[X]$ in each equivalence class under left-composition. Beyond this, figuring out the exact admissible left-factors $\mu$ seems not worth the effort, and we only note that a Siegel function has infinitely many integer values as soon as it has one.

The above classification result allows us to tackle a variety of other classification problems in polynomial arithmetic that were previously inaccessible. Our first application is in arithmetic dynamics. Here, a polynomial iterate $f^{\circ n}(X)$ is called {\it newly reducible} if it is reducible and all previous iterates $f^{\circ j}(X)$, $j<n$ are irreducible. The phenomenon of newly reducible iterates has been studied e.g.\ in \cite{Goksel}, \cite{IJOSS}, and is motivated partly by the much-investigated phenomenon of ``dynamical irreducibility" (also called ``stability", cf.\ \cite{Odoni}) of a polynomial $f$, i.e., simultaneous irreducibility of all iterates of $f$. 
We apply our results to the question of newly reducible iterates ``in fibers" of a given polynomial $f$. For this, we call $f^{\circ n}(X)-t_0$ newly reducible if it is reducible and all $f^{\circ j}(X)-t_0$, $j<n$ are irreducible. This is not a significant generalization of the above notion, since indeed reducibility of $f^{\circ n}(X)-t_0$ is equivalent to reducibility of $\tilde{f}^{\circ n}(X)$ for $\tilde{f}:= (X-t_0)\circ f \circ (X+t_0)$. It is, however, natural to ask how often reducibility phenomena occur as $t_0$ varies. We have the following result, completely classifying which rational polynomials allow a newly reducible $n$-th iterate over infinitely many integer values $t_0$. It extends \cite[Corollary 1.3]{BKN26}, which showed a nonexistence result in absence of composition factors of degree $2$ or $4$.
	
	\begin{theorem}
		\label{cor:newlyred}
		Let $f\in \mathbb{Q}[X]$ be such that there exists an $n\ge 2$ for which $f^{\circ n}(X)-t_0$ is newly reducible for infinitely many $t_0\in \mathbb{Z}$. Then $f$ and $n$ fulfill one of the following, up to conjugating $f$ by linear polynomials over $\mathbb{Q}$:
		\begin{itemize}
			\item[a)] $n=2$, and $f=(aX^s h(X)^2 + b)^2$ for $a\in \mathbb{Q}^\times$, $b\in \{0,1\}$, $s\ge 1$ odd and $h\in \mathbb{Q}[X]$ with $h(0)\ne 0$.
			\item[b)] 
			$n=3$, and $f = (\frac{1}{2}T_m - 1)^2-2$ for $m\ge 1$ odd, where $T_m$ denotes the normalized $m$-th Chebyshev polynomial.
			\item[c)] $n=2$, and $f=(X^2(X+4)^2)\circ (aX^s h(X)^2)$, where $a\in \mathbb{Q}^\times$ is of the form $a=2u^2-v^2$, $s$ is odd and $h\in \mathbb{Q}[X]$ with $h(0)\ne 0$.
			%
			%
			%
			%
			%
			%
			%
			\end{itemize}
			Conversely, if $f$ is of one of the above forms, then some conjugate $\tilde{f}:=\lambda^{-1}\circ f\circ \lambda$, $\lambda\in \mathbb{Q}[X]$ linear, has infinitely many newly reducible $\tilde{f}^{\circ n}(X)-t_0$, $t_0\in \mathbb{Z}$.
		\end{theorem}
	
	Note that the ``base case" $m=1$ of Case b) is conjugate over $\mathbb{Q}$ (via the linear polynomial $4X+2$)  to $X^2-1$. That this polynomial $f$ has infinitely many integer values $t_0$ with $f^{\circ 3}(X)-t_0$ newly reducible is known from \cite{Goksel}.
	
	While the strength of our main results lies in showing that polynomial iterates with infinitely many newly reducible fibers are of a very restricted shape, the converse direction in Theorem \ref{cor:newlyred} also deserves attention. Notably, Theorem \ref{cor:newlyred}b) and c) show that all pairs $(3, 2m)$ with $m$ odd and all pairs $(2,4m)$ with $m$ odd are among the pairs $(n,d)$ such that there exist infinitely many polynomials $f\in \mathbb{Q}[X]$ of degree $d$ with newly reducible $n$-th iterate $f^{\circ n}$. Determination of such values was posed as an open problem in \cite[Question 1.2]{IJOSS}, and the above cases seem to be new.

	Our next application is to simultaneous reducible values of two polynomials. Concretely, given two polynomials $f,h\in \mathbb{Q}[X]$ of degree $>1$, can it happen that $f(X)-a$ and $h(X)-a$ are simultaneously reducible for infinitely many integer values $a\in \mathbb{Z}$? Simultaneous reducibility certainly occurs if $a$ is itself a value of both $f$ and $h$, and therefore this question naturally generalizes the question of which polynomials allow an infinite set of common integer values $f(x)=h(y)$, $x,y\in \mathbb{Q}$. The latter question was fully answered in a landmark paper by Bilu and Tichy \cite{BT}. We extend their result to answer the simultaneous reducibility question.
	\begin{theorem}
		\label{thm:simul}
		Let $f,h\in \mathbb{Q}[X]$ be of degree $>1$. Assume that the set of $a\in \mathbb{Z}$ such that both $f(X)-a$ and $h(X)-a$ are reducible is an infinite set. Then up to replacing $(f,h)$ by $(\lambda\circ f, \lambda\circ h)$ for a linear polynomial $\lambda\in \mathbb{Q}[X]$ and interchanging $f$ and $h$, one of the following holds.
		\begin{itemize}
			\item[1)] (``Trivial case"). $f$ and $h$ have a common indecomposable left-factor $f_1\in \mathbb{Q}[X]$, i.e., $f=f_1\circ f_2$, $h=f_1\circ h_2$ for suitable $f_2,h_2\in \mathbb{Q}[X]$.
			\item[2)] (``Bilu-Tichy case"). $f$ and $h$ have indecomposable left-factors $f_1$ and $h_1$ respectively such that $f_1$ and $h_1$ form a ``standard pair" in the sense of \cite[Section 1.1]{BT} with infinitely many solutions of bounded denominator to $f_1(x)=h_1(y)$.
			\item[3)] $f$ and $h$ have left-factors $f_1$ and $h_1$ respectively such that $(f_1,h_1)$ forms an ``exceptional pair" as defined in Section \ref{sec:pfsimul}. Here, in particular $(\deg(f_1),\deg(h_1))\in\{(2,4), (2,8), (4,4)\}$ in all cases.
			\end{itemize}
	\end{theorem}
	Note that Cases 1), 2) of the above theorem are exactly those in which certain nonlinear left-factors of $f$ and $h$ have an infinite set of common integral values. Case 3) measures the difference between our question and Bilu-Tichy's  problem. That the shape of the polynomials in this case is quite restricted is a testimony to the strength of the classification result Corollary \ref{cor:hilbert_siegel}. We will exhibit examples for each of the families belonging to Case 3, demonstrating that  polynomial pairs in these cases can indeed have an infinite set of simultaneous reducible fibers over integers, while no pair of left-factors of them has infinitely many common integer values.
	
	A special case of Theorem \ref{thm:simul} relates to the problem of ``consecutive" reducible values, and more concretely to the following question: Given an irreducible two-variable polynomial $F(t,X)\in \mathbb{Q}[X]$, how long can a sequence $a+1,\dots, a+k$ of integer values be such that $F(a+j,X)$ is reducible for all $j=1,\dots, k$? While fully effective results seem difficult to establish, an interesting problem is to bound the maximal length of such a reducible sequence for a given polynomial $F$, away from finitely many values. This was achieved in \cite[Proposition 5.5]{DW}, yielding a polynomial upper bound in $\deg_X(F)$ for the maximal $k$ occurring infinitely often as the length of such a sequence. Using our results, we can do much better for the important special case of polynomials of the form $f(X)-t$, corresponding to the question of ``consecutive" reducible fibers under the polynomial map $X\mapsto f(X)$. Indeed, we can show that, apart from one exceptional family, there can only be finitely many occurrences of two consecutive reducible fibers (and as the proof will show, even for the exceptional cases, three consecutive reducible fibers occur only finitely often, cf.\ Remark \ref{rem:gaps}). In particular, for $t$ chosen sufficiently large (with the bound depending on $f$, unfortunately in a non-effective way), one of the polynomials $f(X)-t_0$ and $f(X)-(t_0+1)$ is then always irreducible. More precisely:
	\begin{corollary}
		\label{cor:consec}
		Let $f\in \mathbb{Q}[X]$ be a polynomial which does not have a left-composition factor $\mu\circ D_{2p}(X,\alpha)$ for $\mu\in \mathbb{Q}[X]$ linear, $\alpha\in \mathbb{Q}^\times$ and an odd prime $p$ (where $D_{d}(X,\alpha)$ is the degree-$d$ Dickson polynomial with parameter $\alpha$). Then given any $B> 0$, the set  of values $a\in \mathbb{Z}$ such that $f(X)-a$ is reducible over $\mathbb{Q}$ contains only finitely many pairs of distance $\le B$.
		\end{corollary}
		\begin{remark}
		The exception on left-factors $D_{2p}(X,\alpha)$ is necessary. E.g., $D_{2p}(X,1)-a = T_{2p}-a$ becomes reducible at all values $a$ of $T_2$ and of $T_p$, and moreover the equation $T_2(X)=D_p(Y)-4$ is of the form $X^2 = (T_{p}(Y)-2 =) (Y-2)g(Y)^2$ with $g\in \mathbb{Z}[X]$, and thus  has infinitely many rational solutions with integral $X$-coordinate. Therefore any $f$ with left-factor $T_{2p}$ has infinitely many pairs of reducible integral values of distance $4$. More generally, $D_{2p}(X,\alpha)  = D_2(X,\alpha^p)\circ D_p(X,\alpha) = D_p(X,\alpha^2)\circ D_2(X,\alpha)$ for $p$ odd, and due to $D_p(X,\alpha^2)$ having a finite branch point at $2\alpha^p$, the equation $D_2(X,\alpha^p) = D_p(Y,\alpha^2) - 4\alpha^p$ becomes $X^2 (=D_p(Y,\alpha^2)-2\alpha^p) = (Y-2\alpha)g(Y)^2$ with $g\in \mathbb{Z}[X]$, giving infinitely many pairs of reducible integral values of distance $4\alpha^p$.
		\end{remark}
	
	\subsubsection*{Ideas of proof}
	The proof of Theorem \ref{thm:main} is contained in Section \ref{sec:pfmain}. The analog of the theorem for {\it polynomial} pairs $(f,g)$ was shown (in fact, over arbitrary characteristic zero base fields) in \cite{BKN26}, with an approach based among others on a ``diagonality theorem" (\cite[Theorem 4.1]{BKN26}), originally due to Behajaina and Neftin: if $(f,g)$ is a minimally reducible pair and $f=h\circ f_r$ for an indecomposable polynomial $f_r$, then the kernel $\ker(\textrm{Mon}(f)\to \textrm{Mon}(h))$ of the restriction map on monodromy groups has a socle which is (nearly) diagonal inside the direct product $\textrm{Mon}(f_r)^{\deg(h)}$. The result was then deduced by showing (through in-depth permutation group-theoretical investigations) that, on the contrary, the relevant monodromy groups $\textrm{Mon}(f)$ essentially never fulfill such a diagonality property. It was also remarked (\cite[Remark 6.11(2)]{BKN26}) that for minimally reducible pairs of a  polynomial and a proper Siegel function over $\mathbb{Q}$, there exists a family of examples of degrees $\deg(f)=8$, $\deg(g)=16$ (see also Theorem \ref{thm:8t47} of the present paper) which violates the diagonality property, seemingly implying that the same method cannot be used in this scenario. A decisive new insight is that, surprisingly, this example together with two further ones of the same degree $\deg(f)=8$, actually constitutes the only violations and outside of those the diagonality property continues to hold, see Theorem \ref{thm:weak_diag}. After this, a lot of the machinery of \cite{BKN26} is applicable, although the final steps in the proof require an additional effort for treating certain constellations of right-composition factors of $f$.
	
	The applications reduce, in one way or another, to the investigation of compositional left-factors of polynomials together with their ``partner Siegel functions" (in a minimally reducible pair) if those exist. In particular, apart from the results of Theorem \ref{thm:main}, they make extensive use of Ritt's theorems on polynomial decomposition. Certain parts, such as the proof of Theorem \ref{thm:simul}, also require extensions of Ritt's theorems to Siegel functions, and in particular rely on combinatorial arguments bounding the genus of fiber products of two Galois covers.
	
		It should be made clear that, unlike the special case of odd degree polynomials solved (among others) in \cite[Theorem 1.2]{BKN26}, the result of Corollary \ref{cor:hilbert_siegel} (and hence also all of its further applications) is specific to $\mathbb{Q}$ and does not generalize to other number fields in a straightforward way. This is because the property of Siegel functions over $\mathbb{Q}$ to have poles of the same order is crucially used in the proof of Theorem \ref{thm:main} in Section \ref{sec:pfmain}, and this property does not hold for Siegel functions over arbitrary number fields. The examples of minimally reducible polynomial-Siegel pairs over arbitrary number fields known to the author at this point are considerably more diverse than the situation over $\mathbb{Q}$. At this point, I am aware of solvable examples in polynomial degrees $4$, $6$, $8$, $12$ and $24$, as well as a considerable list of nonsolvable ones. Conjecturally, the final list will be finite and explicit, although giving a concrete conjecture seems premature at this point.

	{\bf Acknowledgement/ AI disclosure}: Some important ingredients for the proofs in Section \ref{sec:pfmain} where obtained through extended discussions with Claude Opus 5.5.
	  In all cases, the author identified the main proof strategy; however, important arguments for the exclusion of $p_r=4$ in Theorem \ref{thm:weak_diag}iii) as well as for Theorem \ref{thm:strong_diag} were first discovered by AI. In all such instances, arguments were independently verified and improved. Claude was also used for proofreading and spell-checking. Magma computations, whenever mentioned, were exclusively carried out by the author.

	\section{Some prerequisites}
	
	In order to avoid doubling up on the exposition of the predecessor paper \cite{BKN26} too much, we recall only the central notions and concepts needed in the following. See \cite[Section 2]{BKN26} for a more in-depth introduction.
	
	\subsection{Monodromy groups}
	The (arithmetic) monodromy group of a rational function $f\in k(X)$ over a field $k\supseteq \mathbb{Q}$ is $\textrm{Mon}_k(f):=\textrm{Gal}(f(X)-t/k(t))$. We will use the shorthand $\textrm{Mon}(f)$ when the base field is clear, and only occasionally switch to the geometric monodromy group $\textrm{Mon}_{\overline{k}}(f)\trianglelefteq \textrm{Mon}_k(f)$. The group $\textrm{Mon}_{\overline{k}}(f)$ has a generating genus-$0$ system (in its transitive action on the roots of $f(X)-t$), consisting of inertia group generators at the branch points of $f$. In particular, if $f$ is a polynomial, resp.\ a proper Siegel function over $\mathbb{Q}$, then $\textrm{Mon}_{\overline{k}}(f)$ contains a cyclic transitive subgroup, resp.\ a cyclic group with two orbits of equal length, namely the inertia group at infinity.
	By the Galois correspondence, it is primitive if and only if $f$ is functionally indecomposable over $k$.  Monodromy groups of indecomposable polynomials and Siegel functions were fully classified by M\"uller in \cite{Mue95} and \cite{Mue}.
	Due to the focus of Theorem \ref{thm:main} on solvable groups, it is useful to recall the (much more elementary) fact that in the solvable case, the geometric monodromy group of an indecomposable polynomial is necessarily cyclic or dihedral of prime degree, or $S_4$. In fact our proof of Theorem \ref{thm:main} in Section \ref{sec:pfmain} does not rely on the above classification results beyond this fact.
	
	\subsection{Siegel functions}
	\label{sec:siegel}
	\begin{definition}
		\label{def:siegel}
		Let $k$ be a number field. A rational function $g\in k(X)$ with at most two poles is called a Siegel function. If $|g^{-1}(\infty)|=2$, then we call $g$ a proper Siegel function. In the special case $k=\mathbb{Q}$, a Siegel function $g\in \mathbb{Q}(X)$ is called a Siegel function over $\mathbb{Q}$ if either $|g^{-1}(\infty)|=1$ or both poles of $g$ have the same order.
		\end{definition}
The relevance of Definition \ref{def:siegel} is due to Siegel's theorem: If $\psi: X\to \mathbb{P}^1_k$ is a cover of irreducible curves over $k$ such that $X$ has infinitely many $k$-rational points $x_0$ with $\psi(x_0)\in \mathcal{O}_k$ integral, then $X\cong \mathbb{P}^1_k$ and $\psi$ is given by a Siegel function. Conversely, if $g\in k(X)$ is a Siegel function, then up to a finite extension of the base field $k$, one obtains infinitely many $x_0\in k$ for which the denominator of $g(x_0)$ has bounded height depending only on $g$ - in particular, there exists a constant $C\in k$ such that $Cg$ has infinitely many integral values. This converse only requires Dirichlet's unit theorem, and the base field extension is never necessary when $|g^{-1}(\infty)|=1$. If $k=\mathbb{Q}$ and one wishes to avoid the base field extension (i.e., obtain infinitely many points of bounded denominator over $\mathbb{Q}$ itself), then necessarily both poles of $g$ are of the same order, and are in fact real and conjugate to each other\footnote{One could thus require the reality of the poles in the definition of ``Siegel over $\mathbb{Q}$". We will not do this in this paper, since our arguments for the most part require only the property of equal pole orders; however, when computing actual examples of Siegel function, this extra property needs to be kept in mind.} - i.e., up to linear transformations, $g(X) = \frac{g_1(X)}{(X^2-c)^d}$ for some $c, d\in \mathbb{N}$, with $c$ not a square. Conversely, this shape guarantees infinitely many rational values of bounded denominator, as easily seen using partial convergents of the continued fraction expansion of $\sqrt{c}$ via the Pell equation.

Using the definition of Siegel functions, it is a standard consequence of the Galois correspondence (for subcovers of a Galois cover of $\mathbb{P}^1_k$) that, given an irreducible two-variable polynomial $F(t,X)\in k[t,X]$, the set of integral values $t_0\in \mathcal{O}_k$ rendering $F(t_0,X)$ reducible is contained in the union of a finite set and the union $\cup_i (\varphi_i(k)\cap \mathcal{O}_k)$ for finitely many Siegel functions $\varphi_i$. Concretely, it suffices to take the finite set indexed by the intransitive subgroups $D$ of the Galois group of $F(t,X)$ over $k(t)$ with the property that the subcover $\Omega^D\to \mathbb{P}^1_k$ of the Galois closure $\Omega\to \mathbb{P}^1_k$ of (the cover of the $t$-line given by) $F(t,X)$ is birational to a cover $x\mapsto \varphi_D(x)$ given by a Siegel function, and one may even restrict to subgroups $D$ which are maximal w.r.t. this intransitivity property (e.g., \cite[Corollary 2.5]{KN24}).
When $F(t,X)=f(X)-t$ for a polynomial $f\in k[X]$, we will denote the finite set of (equivalence classes under right-composition with linears of) such Siegel functions $\varphi_D\in k(X)$ contributing infinitely many reducible values, i.e., such that $\mathcal{R}_f \cap (\varphi_D(k)\cap \mathcal{O}_k)$ is infinite, by $\Psi_{f,k}$ (and omit the index $k$ when $k=\mathbb{Q}$). Since subcovers of the cover of $\mathbb{P}^1_k$ given by a rational function $\varphi$ correspond to compositional left-factors of $\varphi$, the maximality property above means that the $\varphi_D\in \Psi_{f,k}$ have the property that the curve $f(X)-\varphi_D(Y)$ is reducible whereas $f(X)-\varphi_1(Y)$ is irreducible for every $\varphi_1, \varphi_2\in k(X)$ with $\varphi_D=\varphi_1\circ \varphi_2$, $\deg(\varphi_2)>1$.

\subsection{Minimal reducibility and polynomial decompositions}
	\begin{definition}
		\label{def:redpair}
		Let $k$ be a field and $f,g\in k(X)$ rational functions. The pair $(f,g)$ is called a reducible pair if the curve $f(X)=g(Y)$ is reducible over $k$. A reducible pair is called minimally reducible if for any proper left-factor $f_1\in k(X)$ of $f$ (i.e., such that $f=f_1\circ f_2$ with $f_2\in k(X)$ of degree $>1$) and any proper left-factor $g_1$ of $g$, the pairs $(f_1,g)$ and $(f,g_1)$ are irreducible.
		\end{definition}
		
	The minimal reducibility notion relates naturally to the problem of enumerating functional decompositions, and thereby in the case of polynomials in particular to Ritt's theorems, cf.\ \cite{Ritt}, \cite{ZM}.
	These state that any two decompositions of a polynomial $f$ into indecomposables can be transformed into each other by performing finitely many ``Ritt moves", that is, replacing the decomposition $f= (...)\circ a\circ b\circ (...)$ by $f=(...)\circ c\circ d\circ (...)$ for $a\circ b=c\circ d$ where $\gcd(\deg(a), \deg(c))=1=\gcd(\deg(b),\deg(d))$ (Ritt's First Theorem); and moreover any such $(a,b,c,d)$ belongs to one of only two explicitly given types (Ritt's Second Theorem). 
	Since we will repeatedly use Ritt's theorems over the base field $\mathbb{Q}$, we state a version working over arbitrary fields of characteristic $0$, see, e.g., \cite{Gathen}.
	\begin{theorem}
		\label{thm:ritt2}
		Let $k$ be a field of characteristic $0$, and $A,B,C,D\in k[X]$ with $A\circ B = C\circ D$, $\deg(A)=\deg(C)$ and $\gcd(\deg(A),\deg(B))=1$. Then, up to replacing $(A,B,C,D)$ by $(\lambda_1 \circ A\circ \mu_1, \mu_1^{-1}\circ B\circ \lambda_2, \lambda_1\circ C\circ \mu_2, \mu_2^{-1}\circ D\circ \lambda_2)$ for linear $\lambda_i,\mu_j\in k[X]$ and possibly switching the order of $(A,B)$ and $(C,D)$, one of the following holds:
		\begin{itemize}
			\item[1)] $(A,B,C,D) = (D_m(X, \alpha^\ell), D_\ell(X,\alpha), D_\ell(X, \alpha^m), D_m(X,\alpha))$, for integers $\ell, m$, where $D_m(X,\alpha)$ is the degree-$m$ Dickson polynomial with parameter $\alpha\in k^\times$, i.e. the unique polynomial $f$ with $f(X+\alpha/X) = X^m + (\alpha/X)^m$.
			\item[2)] $(A,B,C,D) = (X^n, X^s h(X^n), X^s h(X)^n, X^n)$ for $h\in k[X]$ with $h(0)\ne 0$ and $\gcd(n,s)=1$.
			\end{itemize}
		\end{theorem}
		
		We also emphasize the following well-known result  from the realm of Ritt's theorems, originally due to Engstrom (\cite{Eng}, cf.\ also \cite[Lemma 2.9]{ZM}), which will be used repeatedly.
		\begin{lemma}
			\label{lem:engstr}
			Let $k\supseteq \mathbb{Q}$, and $A,B,C,D\in k[X]$ such that $A\circ B = C\circ D$. If $\deg(C)$ divides $\deg(A)$, then $C$ is a left-factor of $A$.
		\end{lemma} 
	
	If there exist no nontrivial Ritt moves (i.e., with $a,b,c,d$ all nonlinear), then $f$ is said to have a {\it unique decomposition}. Similarly, if there exist no nontrivial Ritt moves involving the right-most indecomposable factor of $f$, then $f$ is {\it right-unique}. Among the many technical tools from predecessor works that remain relevant for the proof of Theorem \ref{thm:main},  let us highlight the following (Part a) is contained in \cite[Proposition 2]{Fried73}, b) and c) are in \cite[Proposition 2.17]{BKN26}): 
	\begin{proposition}
		\label{prop:right_unique}
		If $(f,g)$ is a minimally reducible pair over a number field $k$, then 
		\begin{itemize}
            \item[a)] $f(X)-t$ and $g(X)-t$ have the same Galois closure over $k(t)$.
			\item[b)] A right-factor  $f_2$ of $f$ cannot have geometric monodromy group $\textrm{Mon}_{\overline{k}}(f_2)=C_n$ for $n$ composite or $\textrm{Mon}_{\overline{k}}(f_2)=D_n$ for $n\ne 4$ composite.
			\item[c)] If $f$ is a polynomial, it is necessarily right-unique.
			\end{itemize}
		\end{proposition}
		
		Note that the condition $\textrm{Mon}_{\overline{k}}(f_2)\in \{C_n, D_n\}$ corresponds exactly to $f_1$ being linearly related to a Dickson polynomial (including the case of parameter $\alpha=0$).
		
		We note a useful consequence of these observations on minimally reducible pairs for the sets $\Psi_{f,k}$ as defined above, linking the problems considered in Theorem \ref{thm:main} and in Corollary \ref{cor:hilbert_siegel}.
		\begin{lemma}
			\label{lem:Psi}
			Let $f\in k[X]$. The elements of $\Psi_{f,k}$ are exactly the (equivalence classes under right-composition with linears of) Siegel functions $\varphi\in k(X)$ such that $(f_1,\varphi)$ is a minimally reducible pair\footnote{This explicitly includes the case of trivial pairs $(f_1,f_1)$ for an indecomposable left-factor $f_1$ of $f$.} for some left-factor $f_1$ of $f$ and $\varphi(k)\cap \mathcal{O}_k$ is infinite.
			\end{lemma}
		\begin{proof}
			It follows immediately from the definition of minimal reducibility that, for $\varphi\in \Psi_{f,k}$, there must be some left-factor $f_1$ of $f$ such that $(f_1,\varphi)$ is a minimally reducible pair. Conversely, assume that $(f_1,\varphi)$ is a minimally reducible pair, and $\varphi(k)\cap \mathcal{O}_k$ is infinite. By Proposition \ref{prop:right_unique}a), $\varphi:\mathbb{P}^1_k\to\mathbb{P}^1_k$ is a subcover of the Galois closure of $f:\mathbb{P}^1_k\to\mathbb{P}^1_k$, and to show $\varphi\in \Psi_{f,k}$, we only need to verify the maximality property on the subgroup $D$ associated to $\varphi$ - equivalently, that for every proper left-factor $\varphi_1$ of $\varphi$, the associated subgroup is transitive, i.e., $\varphi_1(X)-f(Y)$ is irreducible. Assume otherwise, then $\varphi_1$ would belong to a minimally reducible pair $(L_1,\varphi_1)$ for some left-factor $L_1$ of $f$, and since $L_1$ and $\varphi_1$ have the same Galois closure (hence the same ramification index at $\infty$), contained in the one of $\varphi$, it follows that $\deg(L_1)|\deg(f_1)$, whence $L_1$ is a left-factor of $f_1$ by Lemma \ref{lem:engstr}. But then $(f_1,\varphi)$ could not have been a minimally reducible pair.
			\end{proof}

	   \subsection{Fiber products and Siegel components}
	   \label{sec:fib}
	Ritt's theorems can be reinterpreted  as a classification of polynomial maps whose fiber product is again given by a polynomial. In our applications, we need to consider the situation in which the value sets of two different Siegel functions $\varphi_1,\varphi_2\in \mathbb{Q}(X)$ have infinite intersection. This translates to existence of a Siegel component on a certain fiber product, which we therefore discuss here.
	For a Siegel function $\varphi$ over $\mathbb{Q}$, write $r(\varphi)\in \{1,2\}$ for its number of poles, and $e(\varphi)$ for the order of (any of) these poles. Given Siegel functions $\varphi_1,\varphi_2$, consider the curve $C: \varphi_1(X) = \varphi_2(Y)$, i.e., the fiber product of $\varphi_1(X)=t$ and $\varphi_2(X)=t$ over the $t$-line. By a {\it Siegel component} of $C$, we mean an irreducible component that is defined over $\mathbb{Q}$, absolutely irreducible and of genus $0$, with either one place over $t=\infty$, or two conjugate real places over $t=\infty$. By definition, the cover of a Siegel component of $C$ over the $t$-line is given by a Siegel function over $\mathbb{Q}$.
	Up to equivalence over $\mathbb{C}$, all possibilities for $\varphi_1, \varphi_2$ were determined independently in \cite{Zieve07} and \cite{Pak} (note that over $\mathbb{C}$, Siegel functions are equivalent to Laurent polynomials), extending the work of Bilu and Tichy \cite{BT} on the case where both $\varphi_1$ and $\varphi_2$ are polynomials. In order to avoid potential intricacies arising from the base field descent, we derive the necessary auxiliary results required for our situation ad-hoc.
	
	\begin{lemma}
		\label{lem:abh}
		Each complex place of $C$ over $t=\infty$ (where $t=\varphi_1(X)=\varphi_2(Y)$) has ramification index $\ell=\textrm{lcm}(e(\varphi_1), e(\varphi_2))$. In particular, an irreducible component of $C$ of degree $D$ over the $t$-line has $D/\ell$ places over $t=\infty$.
	\end{lemma}
	\begin{proof} This is a straightforward application of Abhyankar's lemma.\end{proof}
	
	\begin{lemma}
		\label{lem:dividing_index}
		If $C$ has a Siegel component $C'$, $r(\varphi_2)=2$ and $e(\varphi_1)|e(\varphi_2)$, then $\varphi_2  = \varphi_1\circ \rho$ for a rational function $\rho\in \mathbb{Q}(X)$. 
	\end{lemma}
	\begin{proof}
		By Lemma \ref{lem:abh}, $C'$ is of degree $D\le 2\textrm{lcm}(e(\varphi_1), e(\varphi_2)) = 2e(\varphi_2) = \deg(\varphi_2)$ over the $t$-line. But then $C'$ is of degree $1$ over the $Y$-line, meaning that  $X$  is a rational function $\rho(Y)$ of $Y$ on $C'$. This gives $\varphi_1\circ \rho = \varphi_2$ as claimed.
	\end{proof}
	
	\begin{lemma}
		\label{lem:abs_irr}
		Assume that $\varphi_1$ is indecomposable over $\overline{\mathbb{Q}}$ and that the splitting field of $\varphi_1(X)-t$ is not contained in the one of $\varphi_2(X)-t$ over $\overline{\mathbb{Q}}(t)$. Then $C$ is absolutely irreducible. In particular, the conclusion holds if $\textrm{Mon}_{\overline{\mathbb{Q}}}(\varphi_1)$ is not a quotient of $\textrm{Mon}_{\overline{\mathbb{Q}}}(\varphi_2)$.
	\end{lemma}
	\begin{proof}
		Let $y$ be a root of $\varphi_2(X)-t$, and denote the two splitting fields over $\overline{\mathbb{Q}}(t)$ by $\Omega_{\varphi_1}$ and $\Omega_{\varphi_2}$ respectively. 
		The monodromy group of the cover of $C$ over the $y$-line is $\textrm{Gal}(\Omega_{\varphi_1}(y)/\overline{\mathbb{Q}}(y))\cong \textrm{Gal}(\Omega_{\varphi_1}/\Omega_{\varphi_1}\cap \overline{\mathbb{Q}}(y))$, and contains the non-trivial normal subgroup $N=\textrm{Gal}(\Omega_{\varphi_1}/\Omega_{\varphi_1}\cap \Omega_{\varphi_2})$ of $\textrm{Mon}_{\overline{\mathbb{Q}}}(\varphi_1)$. But by the indecomposability assumption on $\varphi_1$, the latter group is primitive, and a nontrivial normal subgroup of a primitive group is necessarily transitive. This shows that $C$ is absolutely irreducible.
	\end{proof}
	
	\begin{lemma}
		\label{lem:domination}
		Let $\varphi_1$ be a polynomial of degree $d$ and $\varphi_2$ with $r(\varphi_2)=2$, and $\gcd(d,e)=1$ where $e:=e(\varphi_2)$. If $C$ is absolutely irreducible of genus $0$, then one of the following holds:
		\begin{itemize}
			\item[1)]  There is a rational point $t_0\in \mathbb{Q}$ such that all $P\in \varphi_2^{-1}(\{t_0\})$ are of multiplicity divisible by $d$, and $\varphi_1=(\kappa X^d + t_0)\circ \lambda$ for $\kappa \in \mathbb{Q}^\times$ and $\lambda\in \mathbb{Q}[X]$ linear.
			\item[2)] $\textrm{Mon}_{\overline{\mathbb{Q}}}(\varphi_2)$ is dihedral of order $2e$ (in its regular action).
			\end{itemize}
	\end{lemma}
	\begin{proof}
		$C$ has degree $N=2de$ over the $t$-line and, by Lemma \ref{lem:abh}, two places over $\infty$, each of ramification index $de$. Due to the genus-$0$ condition, Riemann-Hurwitz leaves exactly $N$ for the Riemann-Hurwitz contribution at finite ramified places. For $t_0\in \overline{\mathbb{Q}}$ and places $P\in \varphi_1^{-1}(\{t_0\}), Q\in \varphi_2^{-1}(\{t_0\})$, denote the ramification indices over $t_0$ by $e_P$ and $e_Q$, respectively.  The number of places of $C$ over $t_0$ is then $\sum_{P,Q} \gcd(e_P, e_Q)\le d\cdot |\varphi_2^{-1}(\{t_0\})|$, with equality if and only if $\gcd(e_P,e_Q)=e_P$ for all $P,Q$. In particular, the Riemann-Hurwitz contribution at $t_0$ is
		$N-\sum_{P,Q} \gcd(e_P, e_Q) \ge d(2e-|\varphi_2^{-1}(\{t_0\})|)$. Summing over all finite $t_0$, these lower bounds add up to $d$ times the Riemann-Hurwitz contribution over finite branch points of $\varphi_2$, which is $d\cdot 2e = N$. So necessarily equality everywhere, i.e., $\gcd(e_P,e_Q)=e_P$ for all $P,Q$, so $e_P$ divides $\gcd\{e_{Q}: Q\in \varphi_2^{-1}(\{t_0\})\}$. But Riemann-Hurwitz also shows that this $\gcd$ can be $>1$ for at most two finite points $t_0$, and if there are indeed two, then $e_Q=2$ for all places $Q$ over those two points.  The (geometric) monodromy group of $\varphi_2$ is then generated by two involutions, hence a dihedral group, necessarily in its regular action, giving Case 2). If there is only one such $t_0$, it must be the only finite branch point of $\varphi_1$ by the above, hence necessarily rational and totally ramified under $\varphi_1$.  This yields Case 1). 
	\end{proof}

	\section{Reducible value sets: Proof of Corollary \ref{cor:hilbert_siegel} assuming Theorem \ref{thm:main}}
	Using the characterization of the set $\Psi_f$ obtained in Lemma \ref{lem:Psi}, we are reduced to determining the minimally reducible pairs $(f_1,g)$ involving a left-factor $f_1$ of $f$, other than the trivial pairs  $(f_1,f_1)$ with $f_1$ indecomposable.
	This was already achieved in \cite{BKN26} for 1) the case when $g$ is also a polynomial (\cite[Theorem 6.1]{BKN26}), yielding only the possibility $\textrm{Mon}_{\mathbb{Q}}(f_1)=D_4$ when $k=\mathbb{Q}$; 2) the case when $f_1$ has nonsolvable monodromy group (\cite[Lemmas 6.9 and 6.10]{BKN26}), yielding only the family of quintics $f_1$ with $\textrm{Mon}_{\mathbb{Q}}(f_1)=S_5$ known from \cite{DF} in case $k=\mathbb{Q}$;\footnote{The parameter conditions in our footnote to Corollary \ref{cor:hilbert_siegel} are exactly those which make the associated degree-$10$ function field extension rational, so that $g$ becomes Siegel over $\mathbb{Q}$, together with the condition $(x,y)\ne (0,0)$ excluding the degeneration to Dickson or $X^5$.} and 3) the case when $f_1$ is solvable of odd degree, yielding no additional examples (\cite[Theorem 5.7]{BKN26}).
	We only need to combine these with Theorem \ref{thm:main}, together with the explicit computations in $\deg(f)=8$ contained in Section \ref{sec:ex}.
	It then follows directly that the cases with $j\ge 1$ in the situation of  Corollary \ref{cor:hilbert_siegel} are exactly those with a left-composition factor as in Case b), or arising as a composition of three quadratics, as parameterized in the first part of Case c). Note that since in the latter case, $f_1$ factors through a $D_4$-polynomial, it picks up the ``Siegel partner" from that subcase as well, leading to $j\ge 2$. On the other hand {\it every} exceptional (i.e., not equal to a left-factor of $f$) Siegel function contributing to the reducible value set of $f$ must form a nontrivial minimally reducible pair with {\it some} left-factor of $f$. Since any given left-factor has at most one Siegel partner, cf.\ Remark \ref{rem:2to1}, the only way to increase the $j$-values is to construct functions $f$ having more than one of the $f_1$ in this short list as left-factors. This can only happen via Ritt moves among left-factors of $f$, and by Ritt's theorem this forces the involved indecomposables to be of coprime degree. This means that one of the involved left-factors must necessarily be of arithmetic monodromy group $S_5$,\footnote{In particular, such a factor can have a Ritt move only with cyclic polynomials.} and must allow successive Ritt moves with at least two quadratics. One verifies from Ritt's second theorem that this can only occur when this composition of quadratics is actually biquadratic, equivalent to $X^4$ up to composition with linears. We may then assume it to be exactly $X^4$ up to composing the degree-$5$ factor with a suitable linear, and now Ritt's second theorem leaves only $\mu\circ X(X+c)^4\circ X^4 = \mu\circ X^4(X^4+c)^4$, as given.
	
	\begin{remark}
		\label{rem:allpts}
		Since Siegel's theorem is not an effective statement, completely quantifying the (finite) difference between the set $\mathcal{R}_f$ of reducible integer values and the union of value sets of the $\varphi\in \Psi_f$ is usually too much to ask for. There are, however, some lucky group-theoretical coincidences leading, in some special cases, to the conclusion that this difference is in fact the empty set.
	Notably, in the special case $f=X^4$, the Vahlen-Capelli theorem asserts that the reducible values are {\it exactly} the values $a=x^2$ or $a=-4x^4$; group-theoretically, this can be explained by the fact that the two subgroups corresponding to the two functions $\varphi\in \Psi_f$ (namely, maximal block stabilizer of index $2$ inside the arithmetic monodromy group $D_4$; and intransitive subgroup of index $4$ generated by a fixed-point free reflection) are the {\it only} subgroups maximal w.r.t.\ intransitivity. The same phenomenon happens in $S_5$, where the maximal intransitive subgroups are exactly the point stabilizer and $2$-set stabilizer (the latter corresponding to the degree-$10$ Siegel function). Among the exceptional octics in Corollary \ref{cor:hilbert_siegel}, the one in Case c) has an analogous property, as verified by direct calculation in its (arithmetic) monodromy group $C_2\wr C_2\wr C_2$: there are exactly three maximal intransitive subgroups, and these therefore must correspond to the quadratic left-factor and the two exceptional Siegel partners. Therefore, if $f=h$ is such an octic, the value sets of the $\varphi\in \Psi_f$ also cover all of $\mathcal{R}_f$.
	\end{remark}

	We take some time to show that the functions belonging to the set $\Psi_f$ in the situation of Corollary \ref{cor:hilbert_siegel} are in fact all needed to obtain $\mathcal{R}_f$, provided only that they do indeed have infinitely many integral values.
	That such a statement indeed requires a proof can be seen, e.g., via the two functions $\varphi_1(X)=X^2$ and $\varphi_2(X)=\frac{(X-2)^2}{X^2-3}$. Both are quadratic, with the latter clearly not factoring through the former, but nevertheless all of the infinitely many integral values of $\varphi_2$ are squares, essentially due to the fact that $-1$ is not a norm in $\mathbb{Q}(\sqrt{3})$.
	\begin{theorem}
		\label{thm:conv_hilbertsiegel}
		Let $f\in \mathbb{Q}[X]$, and let $\varphi\in \Psi_f$ (i.e., $\varphi(\mathbb{Q})\cap \mathbb{Z}$ is infinite, and $\varphi$ is either an indecomposable left-factor of $f$ or one of the functions $\varphi_i$, $i\le j$). Then the set difference  $(\varphi(\mathbb{Q})\cap \mathbb{Z})\setminus (\cup_{\psi\in \Psi_f\setminus\{\varphi\}} \psi(\mathbb{Q}))$ is infinite.
		\end{theorem}
		\begin{proof}
			Choose any other $\psi\in \Psi_f\setminus\{\varphi\}$. If infinitely many integral values of $\varphi$ are also values of $\psi$, the  fiber product $C:\varphi(X)=\psi(Y)$ has a Siegel component. Since $\varphi$ does not factor through $\psi$, any irreducible component of $C$ is of degree $>1$ over the $X$-line.
			First assume $\varphi$ is a polynomial of degree $d$.  In particular, any genus-$0$ component of $C$ is of degree $D>d$ over the $t$-line ($t=\varphi(X)$). But the number of integral values of height $\le N$ of a degree-$D$ rational function is $O(N^{1/D})$, whereas the set of integral values of a degree-$d$ polynomial map $\varphi\in \mathbb{Q}[X]$, assumed to be infinite, is always of exact order of growth $\Theta(N^{1/d})$. This shows that finitely many $\psi\in \Psi_f\setminus\{\varphi\}$ cannot cover a cofinite subset of the integral points of $\varphi$.
			
			Next assume $\varphi$ is a proper Siegel function, hence necessarily of degree $10=2\cdot 5$ or $16=2\cdot 8$ due to Corollary \ref{cor:hilbert_siegel}, and $\psi\in \Psi_f\setminus\{\varphi\}$. By Lemma \ref{lem:dividing_index}, the pole orders of $\psi$ and $\varphi$ cannot divide each other. But, again from the list of exceptional cases in Corollary \ref{cor:hilbert_siegel}, together with the fact that left-factors $\psi$ of $f$ would necessarily be indecomposable (over $\overline{\mathbb{Q}}$), it follows that at least one of $\psi$ and $\varphi$ is indecomposable, and of splitting field not contained in the splitting field of the other (or else, divisibility of pole orders). Lemma \ref{lem:abs_irr} implies that $C$ is absolutely irreducible. Lemma \ref{lem:abh} implies that, in order for $C$ to be Siegel, $\psi$ needs to be a polynomial, of degree coprime to the pole order $e(\varphi)$ of $\varphi$. Lemma \ref{lem:domination} 
			shows that $\psi$ can only have one finite branch point $t_0$ (without loss of generality, $t_0=0$, so that $\psi = aX^p$ for some $p\ge 2$), and at this point all multiplicities in $\varphi^{-1}(\{0\})$ must be divisible by $p=\deg(\psi)$. For $\deg(\varphi)=16$, this leaves only $p=2$, which however contradicts the coprimeness of pole orders. For $\deg(\varphi)=10$, $p=5$ is impossible since $\varphi$ having ramification index $5$ both at infinity and at $0$ would make the associated degree-$5$ polynomial partner solvable, linearly related to $X^5$. Hence again $p=2$ and all ramification indices at $t=0$ are even, which leaves only one possible genus-$0$ configuration, namely with two finite branch points of ramification indices $2$ and $4$. This case does in fact lead to $C$ being Siegel, and thus needs separate treatment. Since the inertia group generator over $t=0$ has exactly one $2$-cycle (and two $4$-cycles), we may assume $\varphi(\infty)=0$ up to right-composition with linear fractionals. Moreover, it can be seen, either from brute force computation or from the fact that the decomposition group at $t\mapsto \infty$, equal to the normalizer $C_5\rtimes C_4$ of an order $5$-subgroup in (the transitive) $S_5\le S_{10}$, interchanges the two $5$-cycles, that the poles lie in $\mathbb{Q}(\zeta_5)$, and hence, from the fact that they are real and algebraically conjugate, they can be assumed to be $\pm \sqrt{5}$ via normalization. So  $\varphi$ has homogenization $\varphi(X/Y) = c\frac{Y^2\varphi_1(X,Y)^2}{(X^2-5Y^2)^5}$, $c\in \mathbb{Q}\setminus\{0\}$, $\varphi_1\in \mathbb{Q}[X,Y]$ homogeneous. Let $N$ be the l.c.m. of denominators of coefficients of $c\varphi_1$. Then clearly $\varphi(x/y)$ is integral as soon as $x^2-5y^2=\pm 1$ and $N|y^2$. There are infinitely many such solutions (namely, solutions of the Pell equation $x^2-5N^2z^2=\pm 1$) of both signs, since $-1$ is a norm in $\mathbb{Q}(\sqrt{5})$. This yields infinitely many integral values each in two different classes modulo squares, and $\psi(X)=aX^2$ can capture at most one of these two classes. The proof is therefore complete.
			%
			%
			\end{proof}
		
	\section{Simultaneous reducible values: Deriving Theorem \ref{thm:simul} and Corollary \ref{cor:consec} from Corollary \ref{cor:hilbert_siegel}}
	\label{sec:pfsimul}
	To prove Theorem \ref{thm:simul}, we need to define the notion of ``exceptional pair" in the theorem. 
	Recall that for $\alpha\in \mathbb{Q}$, $D_4(X,\alpha)=(X^2-2\alpha)^2-2\alpha^2$ denotes the degree-$4$ Dickson polynomial with parameter $\alpha$.
	\begin{definition}
		\begin{itemize}
			\item
		An exceptional pair of Type 1a is a pair of the form $(-X^2, X^4-\alpha X^2)$, where $\alpha\in \mathbb{Q}$.
		\item 
		An exceptional pair of Type 1b is a pair  of the form $(-X^2, \kappa(X^2-\alpha)^2)$, where $\kappa\in \mathbb{Q}_{>0}$ is a non-square and $\alpha\in \mathbb{Q}^\times$ is of the form $\alpha=u^2-\kappa v^2$, $u,v\in \mathbb{Q}$.
		\item An exceptional pair of Type 1c is a pair of the form $(-X^2, X^2\circ X(X+4)\circ aX^2)$ with $a=2u^2-v^2$, $(u,v)\in \mathbb{Q}^2\setminus \{(0,0)\}$.
		\item An exceptional pair of Type 1d is a pair of the form $(-X^2, X(X-4)(X+s)^2\circ cX^2)$ with $c\in \mathbb{Q}^\times$ and $s=c(a^2-2b^2)-2$, $(a,b)\in \mathbb{Q}^2\setminus\{(0,0)\}$.
		\item An exceptional pair of Type 2 is a pair of the form $(D_4(X,\beta), \rho^2 D_4(X,\alpha) + 2\beta^2 - 2\rho^2\alpha^2)$, where $\alpha,\beta\in \mathbb{Q}$, $\rho\in \mathbb{Q}_{>0}$ is a non-square, $\beta\ne \pm \rho\alpha$, and $\beta$ is of the form $\beta = u^2-\rho v^2 + \rho\alpha$ for $u,v\in \mathbb{Q}$.
		\end{itemize}
		\end{definition}

	We prove a slightly strengthened version of 
	Theorem \ref{thm:simul}. For this, given a polynomial $f\in \mathbb{Q}[X]$, and a positive integer $N$, define 
	\begin{equation}
		\label{eq:bound_denom}
		\mathcal{R}_f^{(N)} := \{a=\frac{b}{c}\in \mathbb{Q}\mid b,c\in \mathbb{Z}, |c|\le N, f(X)-a \text{ is reducible} \}
		\end{equation}
		The technical advantage of this definition is that sets of rationals with bounded denominators are mapped, under any fixed linear polynomial $\lambda\in \mathbb{Q}[X]$, into sets of bounded denominators (with a possibly different bound). Since the sets of reducible values of $f$ and of $\lambda\circ f$ are in one-to-one correspondence via $\lambda$, it suffices to prove the following theorem for one pair $(f,h)$ in a given equivalence class under simultaneous left- (and of course also right-)composition with linears.
	
	\begin{theorem}
		\label{thm:simulstrong}
		Let $f,h\in \mathbb{Q}[X]$ be of degree $>1$ and $N\in \mathbb{N}$. Assume that the set $\mathcal{R}_f^{(N)}\cap \mathcal{R}_h^{(N)}$ is infinite. Then up to replacing $(f,h)$ by $(\lambda\circ f, \lambda\circ h)$ for a linear polynomial $\lambda\in \mathbb{Q}[X]$ and interchanging $f$ and $h$, one of the following holds.
		\begin{itemize}
			\item[1)] (``Trivial case"). $f$ and $h$ have a common indecomposable left-factor $f_1\in \mathbb{Q}[X]$, i.e., $f=f_1\circ f_2$, $h=f_1\circ h_2$ for suitable $f_2,h_2\in \mathbb{Q}[X]$.
			\item[2)] (``Bilu-Tichy case"). $f$ and $h$ have indecomposable left-factors $f_1$ and $h_1$ respectively such that $f_1$ and $h_1$ form a ``standard pair" in the sense of \cite[Section 1.1]{BT} with infinitely many solutions to $f_1(x)=h_1(y)$ of bounded denominator.
			\item[3)] $f$ and $h$ have left-factors $f_1$ and $h_1$ respectively such that $(f_1,h_1)$ is an exceptional pair of one of the types 1a, 1b, 1c, 1d and 2 defined above.
		\end{itemize}
	\end{theorem}
	
	\begin{remark}
		\label{rem:exctypes}
		The exceptional types do not fall into Cases 1 or 2, as demonstrated by the following:
		\begin{itemize}
			\item[a)] In $(f,h)=(-X^2, X^4-\alpha X^2)$, the only indecomposable left-factor of $h$ is $X^2-\alpha X$. Common rational values of $-X^2$ and this left-factor would be rational solutions to the Pythagorean equation $(X-\alpha/2)^2+Y^2=\alpha^2/4$, and there can only be finitely many of bounded denominator. On the other hand, all values of the form $t=-\frac{1}{4}(y^2-\alpha)^2$, $y\in \mathbb{Z}$  are common reducible values for $f(X)-t$ and $h(X)-t$.
			\item[b)] In $(f,h)=(-X^2, \kappa(X^2-\alpha)^2)$, the only indecomposable left-factor of $h$ is $\kappa X^2$, and $\kappa X^2 + Y^2 = 0$ has only the trivial solutions. On the other hand, all values of the form $t=-\frac{\kappa}{4}y^2(y^2-4\alpha)$ are reducible values for $h(X)-t$, so all solutions of $z^2=\frac{\kappa}{4}y^2(y^2-4\alpha)$ give common reducible values.
			This reduces to the conic $\alpha = u^2-\kappa v^2$ assumed to have points.
			\item[c)] In Type 1c, the only indecomposable left-factor of $h=X^2\circ X(X+4) \circ aX^2$ is $X^2$.  At the same time, the infinitude of the set of common reducible values of bounded denominators is evident from Theorem \ref{thm:8t28} (simply choose values of the Siegel function $g$ given in that theorem).
			\item[d)] Similarly, in Type 1d, the infinitude of the set of common reducible values is due to Theorem \ref{thm:8t44}, whereas common values of (the unique) left-factors would be solutions to $-Y^2=X(X-4)(X+s)^2$, and hence via parameter change to $U^2-4 + V^2 = 0$; again, there can only be finitely many of bounded denominator.
			\item[e)] Finally in Type 2, the equation for common values of left-factors would be $X^2-(\rho Y)^2  = 4(\beta^2 - \rho^2\alpha^2)$, with only finitely many solutions of bounded denominator. But values $t=-\frac{1}{4} D_4(x,2\beta) = -\frac{1}{4}(x^2-4\beta)^2 + 2\beta^2$ are reducible values for the first component, and via linear translate, values $t=2\beta^2-(\rho^2/4)(y^2-4\alpha)^2$ are reducible values for the second component. Those occur simultaneously whenever $x^2-\rho y^2=4(\beta-\rho\alpha)$.
			\end{itemize}
	\end{remark}
	
	It may be interesting to note that, regarding the four Type-1 subcases in Remark \ref{rem:exctypes}, there is no infinite set of common values (of bounded denominator) attained by both of the polynomials, but there is still an infinite set of values attained by the first while making the second one reducible. In the last case (``Type 2"), there is no such set either, but an infinite set of values making both reducible (without a linear factor).

				We now deduce Theorem \ref{thm:simulstrong}.
\begin{proof}[Proof of Theorem \ref{thm:simulstrong}]
If there are indecomposable left-factors $f_1$ of $f$ and $h_1$ of $h$ such that $f_1(X)=h_1(Y)$ has infinitely many solutions of bounded denominator, then \cite[Theorem 1.1]{BT} shows $f_1=\phi\circ \tilde{f_1}\circ \mu$, $h_1=\phi\circ \tilde{h_1}\circ\nu$ where $(\tilde{f_1},\tilde{h_1})$ is a standard pair, $\phi\in \mathbb{Q}[X]$ and $\mu,\nu$ are linear. The case $\deg(\phi)>1$ is Case 1) of our theorem, whereas $\deg(\phi)=1$ is Case 2. We may thus assume from now on that $f$ and $h$ have no such left-factors. Then infinite subsets of $\mathcal{R}_f^{(N)}\cap \mathcal{R}_h^{(N)}$ must arise as subsets (of bounded denominator) of the intersection $\psi_1(\mathbb{Q})\cap \psi_2(\mathbb{Q})$ for Siegel functions $\psi_1\in \Psi_f$, $\psi_2\in \Psi_h$ such that at least one of $\psi_1, \psi_2$ is not a left-factor of its associated polynomial. Since infinite subsets of bounded denominator correspond to infinite sets $\mathcal{R}_{\lambda\circ f}\cap \mathcal{R}_{\lambda\circ h}\cap\mathbb{Z}$ of actual integer values for suitable $\lambda = \tilde{N}X$, $\tilde{N}\le N$, we may apply
Corollary \ref{cor:hilbert_siegel} (together with the explicit parameterizations in Section \ref{sec:ex}) to obtain all possibilities for such exceptional $\psi_i$, and we have to go through all possibilities to determine when $\psi_1(X)=\psi_2(Y)$ can have a Siegel component. So from now on, always assume (by symmetry) that $\psi_2\in \Psi_h$ is not a left-factor of $h$.

{\bf Case I}: $\psi_1=:f_1$ is an indecomposable left-factor of $f$.

{\bf Case Ia}: $\psi_2$ is a polynomial. While Bilu-Tichy's list of standard pairs can in principle be consulted to help with the treatment of this case, a direct treatment is more illuminating. This case occurs only if $\psi_2$ is the partner in a minimally reducible pair of a left-factor $h_1$ of $h$ which is a composition of two quadratics. Any such $h_1$ is equivalent, up to left- and right-composition with linears, to $D_4(X,\alpha)$ with $\alpha\in \mathbb{Q}$, and we may therefore assume $h_1=D_4(X,\alpha)$ from now on. Then $\psi_2=-\frac{1}{4} D_4(X,2\alpha)=-\frac{1}{4}(X^2-4\alpha)^2 + 2\alpha^2$. We distinguish two further subcases.
\begin{itemize}
	\item[i)] $\deg(\psi_1)=2$. Via right-composition with linears, we may assume $\psi_1=\rho X^2+t_1$, and we have to look for Siegel components of $C: \rho X^2 = -\frac{1}{4} D_4(Y,2\alpha) -t_1$. Over the $Y$-line, such a component is ramified at all roots of $P(Y):=-\frac{1}{4} D_4(Y,2\alpha) -t_1$ of odd multiplicity. If $P$ has four simple roots, this implies genus $1$, so $P$ must be inseparable. This leaves only $t_1=2\alpha^2$, or $\alpha\ne 0$ and $t_1=-2\alpha^2$. In the first case $P$ equals $-\frac{1}{4}$ times a square, and thus $C$ splits into two components $X=\pm (\dots)/(2\sqrt{-\rho})$ over $\mathbb{Q}(\sqrt{-\rho})$. These can only have infinitely many rational points if they are defined over $\mathbb{Q}$, i.e., when $-\rho$ is a square. Then $\psi_1 = -X^2 + 2\alpha^2$ up to scaling from the right. Via simultaneous left-composition of $f$ and $h$ by $X-2\alpha^2$ this gives the exceptional pair of Type 1a.
	In the second case $t_1=-2\alpha^2$ ($\alpha\ne 0$), $C$ equals $\rho X^2 = -(\frac{Y^2}{2})^2 + 2\alpha Y^2$, which is birational to the conic $U^2 + \rho V^2 = 2\alpha$, with $U=\frac{Y}{2}$. The places over $Y=\infty$ (which is the unique $Y$-value extending $t=\infty$) are real and conjugate if and only if $\sqrt{-\rho}$ is non-rational and real, i.e. when $-\rho>0$ is a non-square; and in that case infinitely many rational points of bounded denominator exist as soon as there is one such point. Up to renaming $-\rho^{-1}=:\kappa$ and composing from the left with $\kappa(X+2\alpha^2)$, this is the exceptional Type 1b.
	
	\item[ii)] $\deg(\psi_1)\ge 3$. As an indecomposable polynomial, $\psi_1$ remains indecomposable over $\overline{\mathbb{Q}}$, and its Galois group is not a quotient of $D_4$. Thus, $C:\psi_1(X)=\psi_2(Y)$ is absolutely irreducible by Lemma \ref{lem:abs_irr}. Assume $C$ is Siegel. If $\deg(\psi_1)$ were even, then every fiber $\psi_1^{-1}(\{t_0\})$ would have to contain (one and hence at least) two elements of odd ramification index, since indeed otherwise $\psi_1$ would factor through a quadratic $cX^2+t_0$. But by Abhyankar, for every finite branch point $t_0$ of $\psi_2$, every odd-index point in $\psi_1^{-1}(\{t_0\})$ becomes a branch point of the cover of $C$ over the $X$-line, whose geometric monodromy group is a subgroup of $\textrm{Mon}_{\overline{\mathbb{Q}}}(\psi_2)$. This either leads to a Siegel function with geometric monodromy group $C_4$ and $\ge 2$ finite branch points, or one with group embedding into $D_4$ and $\ge 4$ finite branch points, two of which are of ramification type $(2,2)$.  Both of these contradict the Riemann-Hurwitz formula. If $\deg(\psi_1)$ is odd, then the cover of $C$ over the $X$-line has one pole, i.e., is given by a polynomial. Since a polynomial cover with group embedding into $D_4$ can have at most two finite branch points, the same reasoning as before gives that the fiber $\psi_1^{-1}(\{t_0\})$ over a branch point of $\psi_2$ can have at most two points of odd ramification index - hence only one since $\deg(\psi_1)$ is odd. This implies in particular $t_0=-2\alpha^2$, so that we obtain $\psi_1(X)+2\alpha^2=c(X-x_1)V(X)^2$ for $c,x_1\in\mathbb{Q}$, $V(X)\in \mathbb{Q}[X]$. Hence $\psi_1(x)=z^2-2\alpha^2$ has the solutions $x=x_1+cw^2$, $z=cwV(x)$, $w\in\mathbb{Z}$. Since $Z^2-2\alpha^2$ is the quadratic left-factor of $h_1=D_4(X,\alpha)$, the indecomposable left-factors $\psi_1$ of $f$ and $Z^2-2\alpha^2$ of $h$ fall under (2). 
	\end{itemize}

{\bf Case Ib}: $\psi_2$ is a proper Siegel function. By Corollary \ref{cor:hilbert_siegel}, it is of degree $2e$ with $e\in \{5,8\}$.
If $C:\psi_1(X)=\psi_2(Y)$ is reducible over $\mathbb{Q}$, then the indecomposable polynomial $\psi_1$ forms a minimally reducible pair with some left-factor of $\psi_2$. Either this left-factor is $\psi_2$ itself, which leaves only the case $\psi_1$ nonsolvable of degree $5$: in this case both $\psi_1$ and $h_1$ are polynomial partners in a reducible pair with the same Siegel function $\psi_2$, and since the degree-$5$ subgroups of the relevant monodromy groups $S_5$, $A_5$ are unique up to conjugation, this forces $\psi_1=h_1$, i.e. Case 1). Or $\psi_1$ is itself a left-factor of $\psi_2$, which can only happen for $e=8$.  In fact, by the computations in Section \ref{sec:ex}, we may assume $\psi_2$ to be one of the functions $g$ in Theorem \ref{thm:8t28} or \ref{thm:8t44} and $\psi_1=-X^2$ (the unique indecomposable left-factor of $\psi_2$). 
These yield exactly the exceptional pairs $(f_1,h_1)$ of Type 1c and 1d respectively.

Otherwise, $C$ is irreducible over $\mathbb{Q}$. If it is Siegel, it must be absolutely irreducible, and with $2\gcd(d,e)$ places over $\infty$ ($d:=\deg(\psi_1)$), forcing $\gcd(d,e)=1$. Since $\psi_2$ is furthermore not dihedral, Lemma \ref{lem:domination} applies as follows.

If $e=8$, then $d$ odd, 
and $\psi_2$ needs to have a fiber with all multiplicities divisible by $d$, contradicting $\deg(\psi_2)=16$.
Otherwise $e=5$. As already noted in the proof of Theorem \ref{thm:conv_hilbertsiegel}, the only cycle type in the primitive action of $S_5$ on $10$ points in which all cycle lengths have gcd $>1$ is $(4.4.2)$. By Lemma \ref{lem:domination}, $\psi_1$ is quadratic, say $\psi_1=\kappa X^2$ and the partner polynomial $h_1$ of $\psi_2$ is of the form $c(Y-y_1)(Y-y_2)^4$ (since ramified of index $4$ at $t_0=0$). But then $\kappa x^2 = c(y-y_1)(y-y2)^4$ has the solutions $y=y_1+(\kappa/c)w^2$, $x=w(y-y_2)^2$, $w\in \mathbb{Z}$, and we are in Case 2) of the theorem.

{\bf Case II}: $\psi_1$ is not a left-factor of $f$ either.

{\bf Case IIa}: Both $\psi_1, \psi_2$ are polynomials. Hence they are of degree $4$, and are the partners in minimally reducible pairs of a left-factor $f_1$, $h_1$ of $f$ and $h$ respectively, which are compositions of two quadratics. After normalizing, we may assume that $\psi_2$ is the partner of  $D_4(X,\beta)$ and $\psi_1$ the partner of $\lambda\circ D_4(X,\alpha)$ for a linear $\lambda(X)=sX+c$. By Bilu-Tichy, $\psi_1=\phi\circ a$ and $\psi_2=\phi\circ b$ for a standard pair $(a,b)$. But no standard pair in \cite{BT} has degrees $(4,4)$, and moreover $\deg(\phi)=4$ would imply $f_1$ and $h_1$ are partners to the {\it same} $D_4$-quartic, in which case they are themselves equal up to right-composition, leading again to Case 1). We are thus left with $\deg(\phi)=2$. Since the partner $-\frac{1}{4}D_4(X,2\alpha)$ of $D_4(X,\alpha)$ has the unique quadratic left-factor $-\frac{1}{4}X^2+2\alpha^2$ (up to right-composition with linears), we get the requirement \begin{equation}
	\label{eq:leftfactors}
	s(-\frac{1}{4}X^2+2\alpha^2) + c = -\frac{1}{4}(rX+w)^2 + 2\beta^2.
	\end{equation} 
	Equating coefficients gives $w=0$, $s=r^2$ and $c=2\beta^2-2r^2\alpha^2$. Moreover, plugging \eqref{eq:leftfactors} into the definition of $\psi_1$ and $\psi_2$, we get $\psi_1 = q\circ (r(X^2-4\alpha))$ and $\psi_2 = q\circ (X^2-4\beta)$, for $q=-\frac{1}{4}X^2+2\beta^2$. Since $q(A)=q(B)$ if and only if $A=\pm B$, the components of $C:\psi_1(X)=\psi_2(Y)$ are the conics
	\begin{equation}
		\label{eq:conic}
		Y^2-\epsilon r X^2 = 4(\beta-\epsilon r\alpha), \ \epsilon =\pm 1.
		\end{equation}
		Set $\rho:=\epsilon r$. If $\rho<0$, there are finitely many points of bounded denominator. If $\rho$ is a square, there are two rational places at infinity. If the right side is $0$, the conic decomposes into a pair of lines, rational only if $\rho$ is a square, and then $\psi_1, \psi_2$ are the same $D_4$ partner (up to right-composition), implying we are once again in Case 1) of the theorem. We are left with the requirement $\rho>0$ nonsquare, right side nonzero, and existence of a rational point. When additionally demanding $\beta\ne -\rho\alpha$, this gives exactly the Type 2 exceptional pair $(f_1,h_1)$. The extra assumption is justified, since if $\beta=-\rho\alpha$, then one may verify directly that the quadratic left-factors of $f_1$ and of $h_1$ coincide up to scaling, leading again to Case 1).

{\bf Case IIb}: $\psi_1$ is a polynomial, whereas $\psi_2$ is a proper Siegel function. Necessarily, $\deg(\psi_1)=4$ and $\deg(\psi_2)=2e$ with $e\in \{5,8\}$. In case $e=5$, Lemma \ref{lem:abs_irr} implies that $C:\psi_1(X)=\psi_2(Y)$ is absolutely irreducible.  Lemma \ref{lem:domination} applies to show that $\psi_2$ has a fiber with all ramification indices divisible by $4$, contradicting $\deg(\psi_2)=10$. Thus $e=8$. By Lemma \ref{lem:dividing_index}, $\psi_1$ is a (decomposable) left-factor of $\psi_2$ of degree $4$. By the computations in Section \ref{sec:ex}, this leaves only $\psi_1=-4X^4$ and $\psi_2$ the Siegel function $g$ in Theorem \ref{thm:8t28}.
But then the corresponding Siegel partner $f_1$ of $\psi_1$ is $X^4$ and the one $h_1$ of $\psi_2$ is the function $f$ in Theorem \ref{thm:8t28}. Both $f_1$ and $h_1$ have a left-factor $X^2$, whence we are again in Case 1).

{\bf Case IIc}: Both $\psi_1$ and $\psi_2$ are proper Siegel functions. From Corollary \ref{cor:hilbert_siegel}, the only possible pole orders for such functions are $5$ and $8$. If $\psi_1$ and $\psi_2$ have the same pole orders, symmetric application of Lemma \ref{lem:dividing_index} yields $\psi_1 = \psi_2\circ \rho$ for a linear fractional $\rho\in \mathbb{Q}(X)$. Hence left-factors $f_1$ of $f$ and $h_1$ of $h$ have the same Siegel partner $\psi_1$, and hence $f_1(X)-t$ and $h_1(X)-t$ have the same Galois closure. 
In the case when $f_1$ is nonsolvable of degree $5$, this immediately forces $f_1=h_1$ since there is only one conjugacy class of index-$5$ subgroups in $S_5$ or $A_5$ respectively. In the case where $f_1$ is an octic as in Section \ref{sec:ex}, Magma computation as in Remark \ref{rem:2to1} shows that in case $f_1\ne h_1$, they nevertheless have a common  left-factor, showing that we are in Case 1 of the theorem after all. 
Finally, if $\deg(f_1)=5$ whereas $\deg(h_1)=8$, note that $\psi_1(X)=\psi_2(Y)$ must be absolutely irreducible by Lemma \ref{lem:abs_irr}, since $\psi_1$ is then indecomposable. But by Lemma \ref{lem:abh}, $\psi_1(X)=\psi_2(Y)$ has $2\cdot 2=4$ places over $t=\infty$, showing that it cannot be Siegel.
\end{proof}

		The following strengthening of Corollary \ref{cor:consec} now becomes a more or less straightforward consequence of Theorem \ref{thm:simulstrong}.
	\begin{theorem}
		\label{thm:consec_strong}
			Let $f\in \mathbb{Q}[X]$ be a polynomial which does not have a left-composition factor $\mu\circ D_{2p}(X,\alpha)$ for $\mu\in \mathbb{Q}[X]$ linear, $\alpha\in \mathbb{Q}\setminus\{0\}$ and an odd prime $p$. Then given any $B> 0$ and $N\in \mathbb{N}$, the set  $\mathcal{R}_f^{(N)}$ contains only finitely many pairs of distance $\le B$.
		\end{theorem}
		
		Before moving to the proof, we include one more elementary lemma.
		\begin{lemma}
			\label{lem:shifted_poly}
			Let $\psi(X)\in \mathbb{Q}[X]$ be a polynomial of degree $n\ge 2$ and $b\in \mathbb{Q}$.
			If the curve $C_b:\psi(X)=\psi(Y)+b$ has a Siegel component, then $b=0$.
		\end{lemma}
		\begin{proof} 
			 The $n$ complex places of $C_b$ over $t\mapsto \infty$ (where $t=\psi(X)$) are then given by the directions $x/y = \zeta$ for $\zeta\in \mu_n$. A component of degree $D$ over the $t$-line must carry $D/n$ of them. If $C$ is a Siegel component of $C_b$, then $D/n\in \{1,2\}$, and in case $D/n=2$, the two poles would have to be real and conjugate. But the only real directions $\zeta$ are $\pm 1$ which are not conjugate to each other. So $D=n$, meaning that the cover of $C$ over the $X$-line is of degree $1$. 
			We have thus shown  that $\psi+b = \psi\circ \rho$ for a rational function $\rho$, which is necessarily of degree $1$. But $\psi\circ \rho$ and $\psi$ have the same branch points whereas the finite sets of branch points of $\psi$ and $\psi+b$ differ by a shift by $b$. The fact that these two are equal forces $b=0$.
		\end{proof}

		\begin{proof}[Proof of Theorem \ref{thm:consec_strong}]
			 Assume that $f,B,N$ are such that $\mathcal{R}_f^{(N)}$ contains infinitely many pairs of non-zero distance $\le B$. Since there are only finitely many differences in $[-B,B]$ of two rational values of denominator $\le N$, there must be some difference $b\in \mathbb{Q}^\times$ (also of bounded denominator only in terms of $f,B$ and $N$) occurring infinitely often. Therefore Theorem \ref{thm:consec_strong} reduces to the special case $h=f+b$ of Theorem \ref{thm:simulstrong} with a constant $b\in \mathbb{Q}^\times$, and it merely remains to check which of the cases of Theorem \ref{thm:simulstrong} can occur with the additional requirement $h=f+b$.
			 
			 Case 1 is obviously impossible, since if $f$ and $f+b$ had a common indecomposable left-factor $f_1$, then $f$ would have the indecomposable left-factors $f_1$ and $f_1-b$ of the same degree, which are not equivalent under right-composition with linears, contradicting Lemma \ref{lem:engstr}. If Case 3 holds then $f$ must have the first component of some exceptional pair as left-factor, while $f\pm b$ has the second component as left-factor. For subtypes a), b) and c) of Type 1 exceptional pairs, the first component is $-X^2$, while the unique (up to right-composition with linears) indecomposable left-factor of the second component is a quadratic with positive leading coefficient. Then $f$ itself would have a quadratic of positive leading coefficient and one of negative leading coefficient as left-factors, which is impossible. For Type 1d), the indecomposable left-factor of the second component is an $S_4$ quartic, and since $f$ having both a $C_2$- and an $S_4$-indecomposable left-factor is again a contradiction to Ritt's first theorem, this case cannot occur either.
			 For Type 2, $f$ would have to have both the quartic left-factors $D_4(X,\beta)$ and $\rho^2 D_4(X,\alpha)+2\beta^2-2\rho^2\alpha^2\pm b$, and since such a quartic left-factor is unique up to right-composition with linears $cX+d$, in particular its leading coefficient is unique up to multiplying by a fourth power. But one of the two factors has leading coefficient $1$, and the other $\rho^2$ with a non-square $\rho$, a contradiction.
			 
			 We are thus left with Case 2 of Theorem \ref{thm:simulstrong}, i.e., the determination of when $(f_1,f_1'-b)$ can be a standard pair in the sense of \cite{BT}, for  indecomposable left-factors $f_1, f_1'$ of $f$. By Lemma \ref{lem:shifted_poly},  $f_1$ and $f_1'$ must be non-equivalent left-factors, and thus have coprime degrees. This leaves only the standard pairs of the first and third kind in \cite[Section 1.1]{BT}. 
			  Here, the third kind is of the form $(D_m(X,a^n), D_n(X,a^m))$ with Dickson polynomials. Recall however that $f_1, f_1'$ are both indecomposable left-factors of the same polynomial, and hence are the two left-factors of a Ritt move. For $m,n\ge 3$, Ritt's second theorem forces $f_1' = D_n(X,a^m)$ to again be a Dickson polynomial, which is only possible for $b=0$ due to the branch locus of these polynomials being symmetric around zero. For $m=2$, the linear shift $X\mapsto X+2a^n$ makes this type a subcase of the first kind, which we are therefore left to consider from now on.
			 The first kind is of the form $(X^m, aX^r p(X)^m)$. By Ritt's second theorem, this means that $f_1' = aX^r p(X)^m+b$ is of the form $X^s q(X)^m$. Putting the Riemann-Hurwitz contributions of the two fibers of $f_1'$ over $0$ and $-b$ into the Riemann-Hurwitz formula already forces $m=2$ and $r=s=1$, and no further branch points. This forces $f_1'$ to be linearly related to a Dickson polynomial of prime degree $p$, and now it is easy to see that the Ritt move occurs if and only if $f_1$ and $f_1'$ are simultaneously linearly related to the left-factors of $D_{2p}(X,\alpha) = D_2(X,\alpha^p)\circ D_{p}(X,\alpha) = D_p(X,\alpha^2)\circ D_2(X,\alpha)$, which is the exceptional case of our theorem.
			 \end{proof}
			 
			 \begin{remark}
			 	\label{rem:gaps}
			 	As the treatment of the exceptional case $f_1=\mu\circ D_{2p}(X,\alpha)$ in the above proof has shown, the gap value $b$ is moreover uniquely determined by the shape of $f_1$: it needs to equal the difference between the two finite branch points of $f_1$. If $f$ has such left-factors $f_1=:f_1(p)$ for different odd primes $p$, this implies Ritt moves and hence equality of the branch point sets of $f_1(p)$ for all those primes $p$.  In particular there is only one infinitely often occurring gap value $b$ for a given $f$, and therefore triples $(a,a+b,a+2b)$ of reducible values cannot occur infinitely often either (since that would imply infinitely many gaps of size $2b$).
			 	\end{remark}

		\section{Newly reducible iterates: Deriving Theorem \ref{cor:newlyred} from Corollary \ref{cor:hilbert_siegel}}

		\begin{proof}[Proof of Theorem \ref{cor:newlyred}]
			An infinite set of newly reducible values for $f^{\circ n}$, $n\ge 2$, can only arise if $f^{\circ n}$ has some left-composition factor $L$ which belongs to a minimally reducible ``polynomial-Siegel" pair $(L,g)$ over $\mathbb{Q}$ and which is not a left-composition factor of $f^{\circ n-1}$. Such $g$ furthermore cannot itself be a left-composition factor of $f^{\circ n}$ (and hence, a fortiori not of $L$), since otherwise, an indecomposable left-factor of $g$ would necessarily (have degree dividing $\deg(f)$ by the uniqueness of the degree sequence of indecomposables due to Ritt's theorems, and thus) be a left-factor of $f$ by Lemma \ref{lem:engstr}, making $f(X)-g(Y)$ reducible. 
			
			Hence $(L,g)$ is a nontrivial minimally reducible pair. By Corollary \ref{cor:hilbert_siegel}, all candidate polynomials $L$ are of degree $4$, $5$ or $8$, and the degree $5$ is not possible either since then $L$ would be indecomposable and would have to be a left-composition factor of $f$, again by Lemma \ref{lem:engstr}.
			 In the same way, $L$ cannot be a composition of a quadratic and an $S_4$ polynomial (in either order), since indecomposables of both these degrees would have to be among the composition factors of $f$, and thus by Lemma \ref{lem:engstr} such a composition would already have to be a left-factor of $f$. 
			 
			 Hence, $L$ is an (arbitrary, since necessarily of arithmetic monodromy group $D_4$) composition of two quadratics or a (specific) composition of three quadratics, respectively. Lemma \ref{lem:engstr} (with $A=f^{\circ 2}$ or $A=f^{\circ 3}$ respectively) shows that $n\le 2$ and $n\le 3$ respectively in those two cases. 

			Assume now first $n=2$. In the case where $L$ has arithmetic monodromy group $D_4$, or equivalently $L$ is a composition of two quadratics, Lemma \ref{lem:engstr} implies that $f$ must be of the form $f=f_1\circ R$ where $f_1$ is quadratic and $\deg(R)$ is odd; indeed, if $4|\deg(f)$, then again $L$ would be a left-factor of $f$. Moreover $R\circ f_1$ must admit a Ritt move $R\circ f_1 = \tilde{f_1}\circ \tilde{R}$ with $\tilde{f_1}$ quadratic (including the trivial Ritt move when $\deg(R)=1$). 
			Via right-composition with linears we may assume $f_1=X^2$ so that Ritt's second theorem forces $R$ to be of the form $R=aX^s h(X)^2+b$ with $a\in \mathbb{Q}^\times, b\in \mathbb{Q}$ and $s$ odd. If $b\ne 0$, one can replace $f$ by its conjugate $\frac{1}{b^2} f(b^2X)$, which is of the shape $(\tilde{a} X^s \tilde{h}(X)^2 + 1)^2$.
			Conversely, the above shape is also sufficient for $f^{\circ 2}$ to acquire a new $D_4$-left factor $f_1\circ \tilde{f_1}$, and hence an infinite set of reducible values of bounded denominator for $f^{\circ n}$ via the partner of $f_1\circ \tilde{f_1}$ in a minimally reducible pair. Conjugating $f$ by a suitable linear $\lambda = B^{-1} X$, $B\in \mathbb{N}$, makes these values integral, and due to Theorem \ref{thm:conv_hilbertsiegel}, infinitely many of these are indeed not values of any other function $\psi\in \Psi_{f^{\circ n}}$ and therefore in particular are newly reducible values, thus showing the converse assertion in Theorem \ref{cor:newlyred}a).
			
			 Assume next that $f^{\circ 2}$ acquires a new left-factor $L$ of shape $(2,2,2)$ as in Corollary \ref{cor:hilbert_siegel}, i.e., up to conjugation, $L=X^2\circ (X(X+4)) \circ (aX^2)$ for some $a\ne 0$, and more precisely $a=2u^2-v^2$, see Theorem \ref{thm:8t28}.
			 By Lemma \ref{lem:engstr}, $\deg(f)\equiv 4\bmod 8$ and $f$ must have the biquadratic left-factor $f_1:=X^2\circ X(X+4)$. Moreover $f$ must be of the form $f=f_1\circ R$, where $R$ allows a Ritt move with $X^2$ to create the additional left-factor $aX^2$, i.e., $R\circ  X^2 = aX^2 \circ \tilde{R}$, or equivalently $(\frac{1}{a} R)\circ  X^2 = X^2\circ \tilde{R}$. Ritt's second theorem shows $\frac{1}{a}R = X^s h(X)^2$.  Conversely, in this case $f^{\circ 2}$ does acquire a new left-factor as in Theorem \ref{thm:8t28}, and thus infinitely many newly reducible values of bounded denominator via its Siegel partner.  This concludes the case $n=2$.
			 
			 It remains to consider the case $n=3$ with $L$ an octic of shape $(2,2,2)$ as above. Then Lemma \ref{lem:engstr} yields that $\deg(f)\equiv 2\bmod 4$ and the maximal length of a chain of quadratic left-factors in $f$ and $f^{\circ 2}$ is $1$ and $2$ respectively.
			 I.e., $f$ itself must necessarily be of the form $f=f_1\circ R$ where $f_1$ is quadratic and $R\circ f_1$ admits a Ritt move.
			 Setting $f_1=X^2$, we get (via the shape of the $(2,2,2)$ polynomial in Theorem \ref{thm:8t28}) the requirement that $R(X^2) \circ R(X^2)$ has a left-factor $X(X+4)\circ (aX^2)$, implying the two Ritt moves $R\circ X^2 = X(X+4)\circ \tilde{R}$ (for a suitable $\tilde{R}$) and moreover $\tilde{R}\circ X(X+4) = (aX^2) \circ \widehat{R}$. In particular, one gets a Ritt move 
			 \begin{equation}
			 	\label{eq:rittmove2}
			 	R\circ (X^2\circ X(X+4)) = (X(X+4)\circ (aX^2))\circ \widehat{R}.
			 	\end{equation}
			 Here $X^2\circ X(X+4) = D_4(X+2,2) + 8$, i.e., there is a Ritt move with left-factor $R(X+8)$ and right-factor the Dickson polynomial $D_4(X,2)$.  From Ritt's theorems, this forces $R(X+8)=\lambda\circ D_m(X,2^4) = \underbrace{\lambda(4^mX)}_{=:\mu} \circ T_m(X/4)$ for $m$ odd and linear polynomials $\lambda,\mu$. To find $\mu$, we may assume $m>1$,\footnote{For $m=1$, i.e. $R$ linear, the conclusion is much easier, e.g. via comparing the quadratic left-factors in \eqref{eq:rittmove2}.} and observe that in the Ritt move \eqref{eq:rittmove2}, the branch points of the right-factors are mapped into the branch point set of the corresponding left-factor, and consequently $R$ and $X(X+4)\circ (aX^2)$ have the same finite branch point set $\{0,-4\}$. Moreover, this Ritt move corresponds to a composition linearly equivalent to an even degree Dickson polynomial, and hence exactly one of the two finite branch points has a fiber with all multiplicities equal to $2$. Looking at $X(X+4)$, this branch point is clearly $-4$, so $R$ must map $0$ to $-4$. On the other hand $R = \mu\circ T_m((X-8)/4)$ has branch point set $\{\mu(2), \mu(-2)\}$. This already leaves only the two options $\mu(X)=X-2$ or $\mu(X)=-X-2$, and the requirement $R(0)=-4$ implies $\mu(-2) = -4$, i.e., finally $\mu=X-2$. 
		We have thus found that, up to conjugation by linears, $f = (T_m((X-8)/4) - 2)^2$. Conjugating by the linear map $4X+8$ yields the shape given in Theorem \ref{cor:newlyred}b).
			 Conversely, we have seen that for $f$ of this shape, $f^{\circ 3}$ does acquire a new octic left-factor of the shape as in Theorem \ref{thm:8t28},\footnote{Indeed, with $a=-2$ modulo squares, as direct computation shows. This is an admissible value in Theorem \ref{thm:8t28}.} and hence infinitely many newly reducible values of bounded denominator via its Siegel partner.
						\end{proof}

		\section{Proof of Theorem \ref{thm:main}}
			\label{sec:pfmain}
		In this section, we apply and sharpen the methods of \cite{BKN26} to prove (the first  implication of) Theorem \ref{thm:main}. Since the ``converse" part will be dealt with independently via explicit computations in Section \ref{sec:ex}, this suffices to prove the theorem in its entirety.
		
		\subsection{Preparatory results}
		We need some results from linear algebra over the finite field $\mathbb{F}_2$.
		
	For the following lemma about elements with two orbits of the same length in affine linear groups,
	cf.\  \cite[Lemma 3.6]{Mue}. 
		\begin{lemma}
			\label{lem:twocycles_affine}
			Assume $V\cong \mathbb{F}_2^a$ acts regularly on a set $B$, and let $\tau\in N_{\textrm{Sym}(B)}(V) = AGL(V)$ consist of two cycles of length $2^{a-1}$. Let $\overline{\tau}$ be the image of $\tau$ in $GL(V)$. Then the following hold.
			\begin{itemize}
				\item[a)] $a\le 4$.
				\item[b)] If $a=4$, then $\overline{\tau}$ is a single unipotent Jordan block.
				\item[c)] If $a=3$, then $\overline{\tau}\ne 1$.
			\end{itemize}
		\end{lemma}
		\begin{proof}
			Embed $AGL(V)$ into $GL_{a+1}(\mathbb{F}_2)$, with $V$ embedded into $\mathbb{F}_2^{a+1}$ as an invariant hyperplane. With this embedding, $\tau$ is unipotent of order $2^{a-1}$, and if $b\le a+1$ is the size of its largest Jordan block, then $a-1 = \lceil \log_2(b)\rceil$. In particular, $a-1\le \lceil\log_2(a+1)\rceil$, which implies
			$a\le 4$.
			 Moreover,  $a=4$ forces $b=a+1=5$, i.e., $\tau$ is a single Jordan block of size $5$ and its restriction to the hyperplane $V$ is necessarily a single block of size $4$. Finally, if $\overline{\tau}=1$, then $\tau$ is a translation, which has order $\le 2$, leading to a contradiction as soon as $a\ge 3$.
		\end{proof}
		
		$\mathbb{F}_2[\langle\overline{\tau}\rangle]$-modules as in Lemma \ref{lem:twocycles_affine} are cyclic $\mathbb{F}_2$-codes. The following narrows down the structure of the particular modules produced by Lemma \ref{lem:twocycles_affine}.
		\begin{lemma}
			\label{lem:cycliccode}
			Let $L\ge 4$ be even, $R:=\mathbb{F}_2[x]/(x^L-1)$, and let $M\le R$ be an $x$-invariant subspace. Write $L=2^eL'$ with $L'$ odd, and let $R_{x+1}:=\ker(x^{2^e}-1)$. Then the following hold.
			\begin{itemize}
				\item[a)] If $\dim M=4$, $4|L$ and $x^{L/4}$ acts on $M$ as a single Jordan block, then $e=2$ and 
				$M=R_{x+1}=\{m\in R\mid x^4m = m\}$.
				\item[b)] If $\dim M=3$ and $x^{L/2}$ acts nontrivially on $M$, then $e=2$ and $M=(x+1)R_{x+1}$.
		    \end{itemize}
			\end{lemma}
			\begin{proof}
				Firstly, $x^L-1 = (x^{L'}-1)^{2^e}$, and $x^{L'}-1$ is squarefree. By Chinese remainder theorem, $R$ is the direct sum of the $R_{\varphi}:= \mathbb{F}_2[x]/(\varphi^{2^e})$, where $\varphi$ runs over the irreducible factors of $x^{L'}-1$. Moreover $M= \oplus_\varphi (M\cap R_{\varphi})$. Here the $R_{\varphi}$ are local PIDs, with submodules $\varphi^s R_\varphi$, $0\le s\le 2^e$, of dimension $(2^e-s)\deg(\varphi)$. In particular, every nonzero submodule of $R_\varphi$ has dimension divisible by $\deg(\varphi)$, and for $\varphi\ne x+1$ the latter is $\ge 2$. For $j\le e$, the element $x^{L/2^j}-1=(x^{L'}-1)^{2^{e-j}}$ projects to $\varphi^{2^{e-j}}u$ on $R_{\varphi}$, with $u$ a unit. Hence on $\varphi^s R_{\varphi}$, it is nilpotent of index $\lceil (2^e-s)/2^{e-j}\rceil$.
				
				For a), take $j=2$. Since $x^{L/4}$ acts as a single block of size $4$, we have $M=\varphi^s R_{\varphi}$ for some irreducible $\varphi$ with $(2^e-s)\deg(\varphi)=4$ and $(4\ge) 2^e-s > 3\cdot 2^{e-2}$. For $e\ge 3$, this yields a contradiction, so $e=2$ and then $s=0$ and $\varphi=x+1$, showing a).
				
				For b), take $j=1$. Since $x^{L/2}$ acts nontrivially, some summand $\varphi^s R_{\varphi}\subseteq M$ has index $\ge 2$, i.e., $2^e-s > 2^{e-1}$, and dimension $(2^e-s)\deg(\varphi)\le 3$. This immediately forces $e\le 2$.  If $e=1$, then $s=0$ and $\deg(\varphi)=1$, so the summand is $R_{x+1}$ of dimension $2$. But then the remaining dimension $2$ cannot come from any other $\psi\ne \varphi$, a contradiction. Hence $e=2$, and thus $s=1$ and $\varphi=x+1$. The equality of dimensions shows $M=(x+1) R_{x+1}$.
				\end{proof}
		
		\subsection{Weak diagonality}
		
		We now develop the main technical tools necessary for the proof of Theorem \ref{thm:main}.
		The following is the analog of \cite[Theorem 4.1]{BKN26} needed for our purpose. It requires the property of Siegel functions over $\mathbb{Q}$ to have poles of the same order, although it does not strictly require the field of definition of the functions to be $\mathbb{Q}$. Therefore, in all of the following, let $k$ be a number field. Monodromy groups, unless stated otherwise, are (arithmetic) monodromy groups over $k$.
		\begin{theorem}
			\label{thm:weak_diag}
			Let $k$ be a number field, $f\in k[X]$  a polynomial with solvable monodromy group, and $g\in k(X)$ a Siegel function with two poles of the same order such that $(f,g)$ is a minimally reducible pair. Let $f=h\circ f_r$ be decomposable with indecomposable right-factor $f_r$, and let $p_r:=\deg(f_r)$ and $K:=\ker(\textrm{Mon}(f)\to \textrm{Mon}(h))$. Then the following hold.
			\begin{itemize}
				\item[i)] If $p_r$ is an odd prime, then $K$ is diagonal, i.e., $|\textrm{soc}(K)| = p_r$.
				\item[ii)] If $p_r=2$, then $K\cong C_2^a$ for some $a\le 4$. If moreover $\deg(f)\ne 8$, then $K\cong C_2$ or $K\cong C_2\times C_2$.
				\item[iii)] If $p_r= 4$, then $\deg(f)=8$.
				\end{itemize}
			\end{theorem}
			
			\begin{corollary}
				\label{cor:tailreduction}
				If $(f,g)$ is a minimally reducible pair where $f\in k[X]$ has solvable monodromy group and $g\in k(X)$ is either a polynomial or a Siegel function with two poles of the same order, then either $\deg(f)\in \{4,8\}$ or $f=h\circ f_{r-2}\circ f_{r-1}\circ f_r$ where $f_{r-2},f_{r-1}, f_r$ are indecomposable of degrees $(\deg(f_{r-2}), \deg(f_{r-1}),\deg(f_r)) \in \{(2,2,2), (2,4,2), (4,2,2)\}$.
				\end{corollary}
				\begin{proof}[Proof of Corollary \ref{cor:tailreduction} assuming Theorem \ref{thm:weak_diag}]
					This is essentially contained in \cite[Section 5]{BKN26}. Namely, the following is shown via combination of Lemmas 5.2, 5.3 and the proof of Theorem 5.1 of that paper: if $f=h\circ f_r\in k[X]$ is right-unique and has no right-factor equivalent to $X^n$, $n$ composite; or $T_n$, $n\ne 4$ composite (as we may assume by Proposition \ref{prop:right_unique}), while also fulfilling the diagonality conclusion of Theorem \ref{thm:weak_diag}i) or ii), then $f$ itself is either a composition of two quadratics (in particular $\deg(f)=4$), or has one of the right-tail shapes given above.\footnote{In particular, the relevant arguments in \cite{BKN26} do not rely on $g$ being a polynomial.} 
					\end{proof}
				
			We now turn to the proof of Theorem \ref{thm:weak_diag}. 	
			Only the case of $p_r\in \{2,4\}$ in Theorem \ref{thm:weak_diag} requires additional efforts, whereas the case $p_r$ odd follows directly from the results of \cite{BKN26}.
			\begin{proof}[Proof of Theorem \ref{thm:weak_diag}i)]
				This is immediate from \cite[Proposition 4.3 and Lemma 4.4]{BKN26} (note only that due to $g$ having two poles of the same multiplicity, an odd degree right-composition factor of $g$ must necessarily be linearly equivalent to a polynomial).
				\end{proof} 
				
			Before turning to the proof in the case $p_r\in \{2,4\}$, we fix some additional notation.
			In the setting of Theorem \ref{thm:weak_diag}, we denote by $\Omega$ the joint splitting field of $f(X)-t$ and $g(X)-t$ over $k(t)$, by $V:=\soc(K)$ the socle of $K$ (which is always an elementary-abelian $2$-group due to $\textrm{Mon}(f_r)\in \{C_2,S_4\}$), and by $\Omega'\subseteq \Omega$ the fixed field of $V$. Moreover, for a root $y$ of $g(X)-t$, write $k(v):=k(y)\cap \Omega'$, and decompose $g = g'\circ w$ with $v=w(y)$, $g', w\in k(X)$.

			The combination of Lemmas \ref{lem:twocycles_affine} and \ref{lem:cycliccode} can be exploited due to the following observation.
				\begin{lemma}
					\label{lem:kernel}
					In the setting of Theorem \ref{thm:weak_diag}, assume $p_r\in \{2,4\}$, so that $V:=\soc(K)$ is of the form $V=C_2^a$. 
					Assume furthermore that $\Omega'(y)=\Omega$. 
					Set $n:=\deg(f)$ and let $\sigma\in \textrm{Mon}(f)$ be an inertia group generator at $t\mapsto \infty$. Let $\overline{\tau}\in GL(V)$ be the automorphism induced on $V$ by conjugation by $\sigma^{n/2^{a-1}}$. Then the pair $(V,\overline{\tau})$ fulfills the assumptions of Lemma \ref{lem:twocycles_affine}.
					\end{lemma}
					\begin{proof}
					Write $g = g'\circ w$ with $v=w(y)$ as defined above.  Let $B:=w^{-1}(v)$ be the subset of the roots of $g(X)-t$ forming the fiber over $v$. Denote by $G_B:=\textrm{Gal}(\Omega/k(v))$ the setwise stabilizer of $B$ and let $\rho: G_B\to \textrm{Sym}(B)$ be the restriction. Since $v\in \Omega' = \Omega^V$, we have $V\le G_B$. Since $\Omega' \cap k(y) = k(v)$, we have $|V| = [\Omega:\Omega'][\Omega'(y):  \Omega']  = [k(y):k(v)] = |B|$. Since $V$ meets the stabilizer of $y$ trivially, this means that $V$ acts regularly on $B$. The element $\sigma^{n/2^{a-1}}$ has order $2^{a-1}$, and since the order of the inertia group at infinity in the Galois closure of $k(v)/k(t)$ divides $\deg(g')=2n/2^a$, $\sigma^{n/2^{a-1}}$ fixes all conjugates of $v$, and in particular lies in $G_B$. Since it has all cycles of length $2^{a-1} = |B|/2$,
						its restriction to $B$ has exactly two cycles, both of length $2^{a-1}$. This shows the assertion.
						\end{proof}
				
				\begin{remark}
					Lemma \ref{lem:kernel} allows us to use the conclusions of Lemma \ref{lem:twocycles_affine}, and hence of Lemma \ref{lem:cycliccode} in the case $p_r\in \{2,4\}$ of Theorem \ref{thm:weak_diag}, {\it assuming} that the technical extra condition $\Omega'(y)=\Omega$ is fulfilled. That this is indeed the case is the content of \cite[Proposition 4.3]{BKN26} in the case $p_r=2$, whereas for $p_r=4$ it will be verified separately before invoking the lemma.
				\end{remark}
		
		The next proposition deals with Conclusion ii) of Theorem \ref{thm:weak_diag}.
		\begin{proposition}
			\label{prop:2_1}
			In the situation of Theorem \ref{thm:weak_diag}, assume that $p_r=2$. Then $a\le 4$. Moreover, if $a\ge 3$, then necessarily $\deg(f)=8$.
			\end{proposition}
			\begin{proof}
				Let $\Omega$, $\Omega'$, $y$ and $v$ be as introduced above, so that $\textrm{Gal}(\Omega/\Omega')\cong C_2^a$. 
				By \cite[Proposition 4.3(1)]{BKN26}, we have $\Omega=\Omega'(y)$. 
				Lemma \ref{lem:kernel} and Lemma \ref{lem:twocycles_affine} give $a\le 4$.
				
				Assume from now on that $a\ge 3$.
								 Let $X = \{x_1,\dots, x_N\}$ be the roots of $h(X)-t$, so that $K$ embeds naturally into the direct product $C_2^X$. Let $\sigma$ be an inertia group generator at $t\mapsto \infty$ and $\overline{\sigma}$ its image on $X$, assumed to fulfill $\overline{\sigma}(x_i)=x_{i+1}$ via suitable numbering. Let $e_i\in C_2^X$ be nontrivial exactly on the fiber over $x_i$. Then $\sigma e_i\sigma^{-1}=e_{i+1}$, so $e_i\mapsto x^i$ identifies $C_2^X$ with  $\mathbb{F}_2[x]/(x^N-1)$ and conjugation by $\sigma$ with multiplication by $x$. Since $K\trianglelefteq G$, it is an $x$-invariant subspace, and conjugation by $\tau=\sigma^{\deg(f)/2^{a-1}}=\sigma^{2N/2^{a-1}}$ is multiplication by $x^{N/2^{a-2}}$ (recall that $\deg(f)/2^{a-1} = \deg(g)/2^a = \deg(g')$ is an integer, where $g=g'\circ w$ is defined via $v=w(y)$).
				
				We next claim that $N=4N'$ with an odd integer $N'$. This is true since, due to Lemma \ref{lem:kernel}, the conclusion of Lemma \ref{lem:cycliccode} holds (with $L=N$ and $M=K$).
				
				Moreover that lemma shows that, in case $a=4$, $K$ consists of the vectors stable under shift by $4$ positions, i.e., of the form $\sum_i c_i e_i$ where $c_i=c_{(i+4) \bmod N}$ for all $i$. In case $a=3$, $K$ consists of the subset of those vectors fulfilling additionally $c_0+c_1+c_2+c_3=0$ (since indeed, $(x+1)R_{x+1}$ is the image of $1+x$ on the $x^4$-invariant vectors). We use this shape to construct a block system for $\textrm{Mon}(h)$. Call $x_i$ and $x_{i'}$ equivalent if every element of $K$ has the same component entry at both. The equivalence classes are the residue classes $X_r:=\{x_i: i\equiv r \bmod 4\}$; this is obvious for $a=4$, whereas for the shape obtained in case $a=3$, merely note that e.g. the vector $(c_0,c_1,c_2,c_3)=(1,1,0,0)$ separates $i=0$ from $i=2$ and $i=3$, etc. Since $g\in \textrm{Mon}(f)$ acts on $K$ via its image $\overline{g}\in \textrm{Mon}(h)$, $\{X_0,\dots, X_3\}$ is a block system of $\textrm{Mon}(h)$, so $h=h_L\circ h_R$ with $\deg(h_L)=4$ and $\deg(h_R)=N'$.
				
				Suppose that $N'>1$, and choose a complete decomposition refining $f=h_L\circ h_R\circ f_r$. We may assume the right-most factor of $h_R$ to be $f_{r-1}$, of odd prime degree $p$ since $N'$ is odd. Since elements of $K$ are constant on fibers of $h_R$, and hence in particular on fibers of $f_{r-1}$, the image $\overline{K}$ of $K$ in $\textrm{Mon}(f_{r-1}\circ f_r)$ lies in the diagonal of the kernel $C_2^p$. On the other hand, let $N_0$ be the normal closure of $\langle\sigma \rangle\cap \ker(\textrm{Mon}(f)\to \textrm{Mon}(f_1\circ\dots\circ f_{r-2}))$, and let $\overline{N_0}\le \textrm{Mon}(f_{r-1}\circ f_r)$ be the image of $N_0$ under the same projection. This is a normal subgroup containing a $2p$-cycle, and since $f_{r-1}$ and $f_r$ do not form a Ritt move, it follows from \cite[Corollary 3.2(2)]{BKN26} that $\overline{N_0}$ contains all of $C_2^p=\ker(\textrm{Mon}(f_{r-1}\circ f_r)\to \textrm{Mon}(f_{r-1}))$. Hence $\overline{N_0} = C_2\wr C_p$ and $\overline{N_0}/(\overline{K}\cap \overline{N_0})$ is a quotient of a subdirect power of $C_p$ (this is a special case of Lemma \ref{lem:absorption} below). 
				In particular,  $\overline{N_0}/(\overline{K}\cap \overline{N_0})$ is of odd order, and therefore $C_2^p\le \overline{K}$, contradicting the fact that $\overline{K}$ is diagonal. Hence $N'=1$ and $\deg(f)=8$.
				\end{proof}
				
				Finally, we move towards Assertion iii) of Theorem \ref{thm:weak_diag}.
				\begin{lemma}
					\label{lem:pr4}
					In the situation of Theorem \ref{thm:weak_diag}, assume $p_r=4$, $r\ge 2$. Then the following hold:
					\begin{itemize}
						\item[a)] $\Omega'(y)=\Omega$.
						\item[b)] $\soc(K) \cong C_2^a$ with $a\le 4$.
						\end{itemize}
					\end{lemma}
					\begin{proof}
					We have $\soc(K) = K\cap V_4^\mathcal{X}$ where $\mathcal{X}$ is the set of roots of $h$. 
					Let $u=f_r(x)$ for a root $x$ of $f(X)-t$. Since the minimal polynomial $h(X)-t$ of $u$ splits in the fixed field of $K\supseteq \soc(K)$, we have $u\in \Omega'$, whereas $x\notin \Omega'$. Since $k(x)/k(u)$ has no intermediate fields and $\Omega'/k(u)$ is Galois, $\Omega'\cap k(x)=k(u)$. Therefore $\Omega'(x)/\Omega'$ is Galois with group the regular $V_4$, and contains all four roots of $f_r(X)-u$. 
					
					Next, let $k(v)=\Omega'\cap k(y)$. By minimal reducibility, $k(x), k(v)$ and $k(u), k(y)$ are linearly disjoint over $k(t)$. Hence $k(x,v)/k(u,v)$ has degree $4$ and is not linearly disjoint from $k(u,y)/k(u,v)$. Let $M_x\subseteq \Omega'(x)$ be the Galois closure of $k(x,v)/k(u,v)$, and $R:=\textrm{Gal}(M_x k(u,y)/k(u,y))$. Then $R$ is intransitive on the four roots. There are thus only two possibilities: 
					
					Case 1) $R$ may fix a root $x_0$. Then $x_0\in k(u,y)$, so $\Omega'(x)\subseteq \Omega'(y)$. As $G_y:=\textrm{Gal}(\Omega/k(y))$ is transitive on the roots of $h$, we have $\Omega'(y)=\Omega$, and thus Lemma \ref{lem:kernel} becomes applicable. The assertion follows from Lemma \ref{lem:twocycles_affine}. 
					
					Case 2) The only other option is that $R$ has two orbits $\{x_1,x_2\}$ and $\{x_3,x_4\}$. We show that this is impossible. Note that in this case $R$ fixes a root $\rho_0$ of the cubic resolvent  of $f_r(X)-u$. Since $\rho_0\in M_x$, we get $\rho_0\in M_x^R=M_x\cap k(u,y)$. On the other hand $G/\soc(K)$ acts faithfully on the total of $3N$ resolvent roots, and $K':=K/\soc(K)\le S_3^N$. Here $\soc(K')=K'\cap C_3^N$ projects onto $C_3$ on each block since $K$ projects onto $S_4$ on each block. Let $\Omega'' = (\Omega')^{\soc(K')}$. Since $\Omega'/\Omega''$ is an abelian extension, $\Omega''(\rho_0)$ contains all three resolvent roots over $u$.
					In particular, since we may also assume that Case 2 holds with $u$ replaced by any other root $\tilde{u}$ of $h(X)-t$ (or else, the assertion follows from Case 1 using that root), we see that all resolvent roots of $f_r(X)-\tilde{u}$ lie in $\Omega''(y)$. Since the set of all these resolvent roots, when running over the $\tilde{u}$, generates $\Omega'$, it follows that $\Omega'\subseteq \Omega''(y)$.
					 Set $V=\soc(K)$ and $\Sigma=\textrm{Gal}(\Omega/\Omega'')$, then this implies $\Sigma\cap G_y\le V$, where $G_y$ is the stabilizer of $y$. With $k(v'):=\Omega'\cap k(y)$ (which is the fixed field of $VG_y$) and $k(v''):=\Omega''\cap k(y)$, Dedekind's modular law gives $\Omega''(v') = \Omega'$. So $|\soc(K')| = [k(v'):k(v'')]=:3^b$, and the composition factor of $g$ corresponding to $k(v')/k(v'')$ has monodromy group containing a regular normal subgroup $C_3^b$. But also, being an odd degree composition factor of a Siegel function with two poles of the same order, this monodromy group has a transitive cyclic subgroup. It follows (see \cite[Lemma 2.2]{BKN26}) that $b=1$. Consequently, $\soc(K')$ is a diagonal $C_3$ and $K'$ is a diagonal $S_3$. In particular, the image of $K$ under projection to $\textrm{Mon}(f_{r-1}\circ f_r) (\le S_4\wr \textrm{Mon}(f_{r-1}))$ is a normal subgroup with diagonal $S_3$-part. But a diagonal $S_3$ is self-normalizing in the whole direct product $S_3^{\deg(f_{r-1})}$. So $\ker(\textrm{Mon}(f_{r-1}\circ f_r)\to \textrm{Mon}(f_{r-1}))$ must also have diagonal $S_3$-part. Since $f_{r-1}$ and $f_r$ cannot form a Ritt move, this contradicts \cite[Theorems 2.3, 2.4]{BKN26b}, which shows that this $3$-part is always of rank at least $2$. 
					Thus, Case 2) cannot happen and the proof is complete.
						\end{proof}
						
			\begin{corollary}
				\label{cor:pr-1}
				If $p_r=4$ and $r\ge 2$, then $p_{r-1}=2$ and $a=4$.
				\end{corollary}
				\begin{proof}
					The image of $\textrm{soc}(K)$ in a $2$-step block is a nontrivial normal subgroup of $\textrm{Mon}(f_{r-1}\circ f_r)$ inside $V_4^{p_{r-1}}$. But $V_4^{p_{r-1}}$ is minimal normal in this $2$-step monodromy group, see \cite[Corollary 3.2(3)]{BKN26} (this Corollary only requires that there is no Ritt move between $f_{r-1}$ and $f_r$, which is automatic due to the right-uniqueness of a polynomial $f$ involved in a minimally reducible pair, cf.\ Proposition \ref{prop:right_unique}). Hence the above image is all of $V_4^{p_{r-1}}$, of order $2^{2p_{r-1}}$, forcing $p_{r-1}=2$ and $a=4$ due to Lemma \ref{lem:pr4}.
					\end{proof}
					
					The proof of Theorem \ref{thm:weak_diag} is completed via the following
			\begin{proposition}
				If $p_r=4$, then $r=2$ and $\deg(f)=8$.
				\end{proposition}
				\begin{proof}
					By the previous corollary, $a=4$.
					Let $\mathcal{X}$ be the set of roots of $h$, labeled as $\mathcal{X}=\{x_i\mid i\in \mathbb{Z}/N\}$. Denote by $F_i$ the fiber of $f_r$ over $x_i$, and by $V_4^{(i)}$ the Klein group acting on $F_i$. Since $n=\deg(f)=4N$, the element $\sigma^N$ stabilizes each $F_i$ and acts on it as a $4$-cycle. Conjugation by $\sigma^N$ thus induces a transposition $\theta$ on the component $V_4^{(0)}$. Let $e,\theta(e)$ be the points in the support of $\theta$. Moreover $\sigma^{2N}$ lies in $V_4^\mathcal{X}$.
					
					For $j\in \mathbb{Z}/2N$, put $v_j:=\sigma^j e \sigma^{-j}\in V_4^{(j\bmod N)}$. Then $\{v_j, v_{j+N}\}$ is a basis of $V_4^{(j\bmod N)}$, so $\{v_j: j\in \mathbb{Z}/2N\}$ is a basis of $V_4^\mathcal{X}$. Conjugation by $\sigma$ maps $v_j$ to $v_{j+1}$, giving $V_4^\mathcal{X} \cong \mathbb{F}_2[x]/(x^{2N}-1)$, and $\tau:=\sigma^{n/8}$ acts as $x^{2N/4}$.  Via Lemma \ref{lem:kernel}, we may invoke Lemma \ref{lem:cycliccode}a) with $L=2N$, which gives $N\equiv 2\bmod 4$, and $\soc(K) = \{\sum_j c_jv_j \mid c_j =c_{(j+4)\bmod 2N} \text{ for all } j\}$.
						
						For a lift $j$ of $i\in \mathbb{Z}/N$ to $\mathbb{Z}/2N$, the basis vectors on $x_i$ are $v_j, v_{j+N}$, and since $N\equiv 2\bmod 4$, we have $j+N\equiv j+2\bmod 4$. Since  $\soc(K)$ consists of exactly the vectors invariant under shift by four positions, the minimal nonempty supports (as subsets of $\mathbb{Z}/N$) of elements of $\soc(K) \subseteq V_4^N$ are exactly the two $\langle\sigma^2\rangle$-orbits $\mathcal{X}_0, \mathcal{X}_1$, each of size $N/2$. Since supports of minimal size must be permuted under conjugation in $\textrm{Mon}(h)$, $\{\mathcal{X}_0, \mathcal{X}_1\}$ is a block system of size $2$ for that group, and hence $h=h_L\circ h_R$ with $\deg(h_L)=2$ and $\deg(h_R)=N/2$ odd. If $N/2>1$, then plugging this into $f=h\circ f_r$ gives a decomposition of $f$ with $p_{r-1}$ odd, contradicting Corollary \ref{cor:pr-1}. Hence $N=2$ and $\deg(f)=8$.
					\end{proof}
				
		\subsection{Strong diagonality}
		\label{sec:strong}
		Due to Corollary \ref{cor:tailreduction}, we are reduced to excluding certain patterns of three indecomposable right-factors for $f$. In order for this not to escalate into consideration of longer and longer right-tails, we sharpen the diagonality conclusion of Theorem \ref{thm:weak_diag} further to tighten the conditions that come with the case $K\cong C_2\times C_2$.
		
		Throughout this section, we assume $n=\deg(f)\ne 8$, $f=h\circ f_r$ with $p_r:=\deg(f_r)=2$, and $h=h'\circ f_{r-(s-1)}\circ \dots \circ f_{r-1}$ for some $s\ge 2$. Furthermore denote the ``right-tail" $f_{r-(s-1)}\circ\dots\circ f_r$ of $f$ by $F$ and the right-tail $f_{r-(s-1)}\circ\dots\circ f_{r-1}$ of $h$ by $F'$ (so that $h=h'\circ F'$). We index the roots of $f$ by $\mathbb{Z}/n$, with the inertia group generator $\sigma$ at $\infty$ acting as the cyclic shift $i\mapsto i+1$. Let $e=\deg(h')$ and $d=\deg(F)$ (so that $de=n$). Write $L=\ker(G=\textrm{Mon}(f)\to \textrm{Mon}(h'))$, so that $\langle\sigma \rangle\cap L = \langle\sigma^e\rangle$, and let $N_0\le L$ be the normal closure of $\sigma^e$ in $G$. Fix a block $\Delta$ of the partition induced by the composition $f=h'\circ F$, let $G_{\Delta}\supseteq L$ be its stabilizer, and $\pi:G_{\Delta}\to \textrm{Sym}(\Delta)$ the restriction (so that $\pi(G_{\Delta})=\textrm{Mon}(F)$). Note that $\pi(\sigma^e)$ is a $d$-cycle, contained in the transitive normal subgroup $G_0:=\pi(N_0)$ of $\textrm{Mon}(F)$. Let $B:=\ker(\textrm{Mon}(F)\to\textrm{Mon}(F'))$ be the minimal block kernel in $\textrm{Mon}(F)$, and finally set $M:=\pi(N_0\cap K)$ and $Q:=G_0/(G_0\cap B)$.
		
The following is then a special case of \cite[Corollary 2.20]{BKN26}.
\begin{lemma}
	\label{lem:absorption}
	With the above notation, $N_0/(N_0\cap K)$ is a subdirect product of $e$ copies of $Q$, and $G_0/M$ is isomorphic to a quotient of $N_0/(N_0\cap K)$. In particular, the nilpotency class and the exponent of the $2$-Sylow group of $G_0/M$ are upper-bounded by the ones of the $2$-Sylow group of $Q$.
	\end{lemma}
		
		\begin{theorem}
			\label{thm:strong_diag}
			With the notation of this section, assume $p_r=2$ and $K\cong C_2\times C_2$ in Theorem \ref{thm:weak_diag}. Let $\mathcal{X}$ be the set of roots of $h$, indexed by $\mathbb{Z}/(n/2)$, and let  $z=\sigma^{n/2}$ be the power of the inertia group generator $\sigma$ lying in $K$ (in particular, $\textrm{supp}(z) = \mathcal{X}$), and let $k_1$ and $zk_1$ be the two further non-identity elements of $K$. Write $S=\textrm{supp}(k_1)\subseteq \mathcal{X}$, so that $\textrm{supp}(zk_1)=\mathcal{X}\setminus S$. Then the following hold:
			\begin{itemize}
				\item[i)] $\{S, \mathcal{X}\setminus S\}$ is a block system of $\textrm{Mon}(h)$. In particular, $h=q\circ u$ with polynomials $q,u\in k[X]$ such that $\deg(q)=2$.
				\item[ii)] If $|\pi(K)|=4$, then $e=\deg(h')$ is odd.
				\item[iii)] In the situation of ii), let $\Delta'\subseteq \mathcal{X}$ be the image of $\Delta$, so that $|\Delta'|=d/2$. Then $S\cap \Delta'$ is a block of size $d/4$. In particular, $h$ has a right-factor $b$ of degree $d/4$, and the right-factor $F'$ of $h$ attached to $\Delta'$ satisfies $F'=a\circ b$ where $\deg(a)=2$.
				\item[iv)] If moreover $e>1$ in the situation of ii), then $f$ has a complete decomposition $f=q\circ (\cdots) \circ v\circ b\circ f_r$ where $\deg(v)$ is an odd prime. In other words, the left-factor $a$ of $F'$ is involved in a Ritt move with a right-factor of $h'$ (necessarily of odd degree).
				\end{itemize}
			\end{theorem}
			\begin{proof}
			 $K$ is normal and the diagonal element $z$ is central in $G$, being the only element of $K$ of full support, whence $\{k_1, zk_1\}$ is permuted by the conjugation action of $G$. Thus $G$ also permutes the supports of these two elements, showing i). Since the cyclic shift $i\mapsto i+1$ on $\mathbb{Z}/(n/2)$ is contained in $\textrm{Mon}(h)$, we may assume $S=2\mathbb{Z}/(n/2)$. The block $\Delta'$ containing $x_0\in S$ is then the subgroup of order $d/2$, i.e., equals $e\mathbb{Z}/(n/2)$. Now $|\pi(K)|=4$ exactly when the restriction ${k_1}_{|\Delta}$ of $k_1$ to $\Delta$ is different from both $1$ and $z_{|\Delta}$, i.e., when $\Delta'$ meets both $S$ and $\mathcal{X}\setminus S$, i.e., when $e\mathbb{Z}\not{\subseteq} 2\mathbb{Z}$, or equivalently, when $e$ is odd. This shows ii). In this situation, one moreover obtains that $S\cap \Delta' = 2e\mathbb{Z}/(n/2)$ has order $d/4$ and index $2$ in $\Delta'$, showing iii). Finally, to see iv), note that the block systems of $S\cap \Delta'$ and of $S$ give right-factors $b$ and $u=u_1\circ b$ of $h$, with $\deg(u_1)=[S:(S\cap \Delta')]=e$ odd. If $e>1$, then an indecomposable right-factor $v$ of $u$ must be of odd prime degree, and comparing sizes of blocks gives the decomposition in the assertion.
				\end{proof}

		\subsection{Proof of the main theorem}
		We divide the remaining proof of Theorem \ref{thm:main} into several steps. First we show:
		\begin{lemma}
			\label{lem:red_to_4step}
			In the situation of Theorem \ref{thm:main}, assume $n=\deg(f)\ne 8$. Then $f$ has a decomposition $f=h\circ f_{r-3}\circ f_{r-2}\circ f_{r-1}\circ f_r$ where the four right-factors are indecomposable of degrees $(\deg(f_{r-3}),\dots, \deg(f_{r}))=(p,2,2,2)$ for $p$ a prime or $p=4$.
			\end{lemma}
		\begin{proof}
			Recall that, due to Corollary \ref{cor:tailreduction}, we are reduced to considering $f=h\circ f_{r-2}\circ f_{r-1}\circ f_r$ with $(\deg(f_{r-2}), \deg(f_{r-1}), \deg(f_r))\in \{(2,2,2), (2,4,2), (4,2,2)\}$.\footnote{Or indeed $\deg(f)=4$, although it is obvious that this case cannot occur for a proper Siegel function with two poles of the same order, since indeed the only index-$8$ subgroups of $S_4$ and $D_4$ are contained in their point stabilizers.} We may exclude the case $\deg(f)=8$, and since it can directly be verified with Magma that there are no minimally reducible pairs in degree $16$ of the respective shapes, we may assume $\deg(h)>1$. If $|\pi(K)|=4$, then Theorem \ref{thm:strong_diag} would yield a right-tail of $f$ of the degree sequence $(p,2,2)$ or $(p,4,2)$, which however is already excluded. We may therefore assume $|\pi(K)|=2$, and thus search for transitive groups $A$ of degree $8$ or $16$ with all of the following properties.
			\begin{itemize}
				\item[1)] $A$ is the (arithmetic) monodromy group of a polynomial $f_{r-2}\circ f_{r-1}\circ f_r$ of shape as above.
				\item[2)] $A$ does not embed into the symmetric normalizer of a cyclic transitive subgroup, and a maximal block stabilizer of $A$ does not act as $C_4$ on the respective block.
				\item[3)] There is a normal subgroup $G_0\trianglelefteq A$, contained in the {\it geometric} monodromy group of $f_{r-2}\circ f_{r-1}\circ f_r$, such that $G_0$ contains a cyclic transitive subgroup and a central subgroup $M$ of order $2$ such that $G_0/M$ embeds as a quotient into a subdirect product of $\textrm{Mon}(f_{r-2}\circ f_{r-1})$.
				\end{itemize}
			
			The only groups meeting all the requirements are the groups $A$ of transitive group label 8T26 (of shape $(2,2,2)$, necessarily with geometric monodromy group 8T17), 16T773 or 16T1060 (of shape $(4,2,2)$, both with geometric monodromy group 16T773). The latter two however have image $GL(2,3)$ in the action on a minimal block system, and since this is not a monodromy group of a polynomial defined over $\mathbb{Q}$ (e.g., \cite[Appendix A]{BKN26b}), it may be excluded for our purposes.\footnote{This is the only place where we have in fact invoked the field of definition being $\mathbb{Q}$, rather than merely the (a priori weaker) condition that the Siegel function has two poles of the same order. Presumably, this case could also be dealt with via suitable analogs of the following three lemmas for tail shape $(\star, 4, 2, 2)$, but since permutation degrees become quite large for carrying out the necessary Magma checks, we content ourselves with invoking $k=\mathbb{Q}$ for once.}
			We are therefore reduced to $f=h\circ f_{r-2}\circ f_{r-1}\circ f_r$ where $\deg(h)>1$ and all three right-factors are quadratic; in fact $A=\textrm{Mon}(f_{r-2}\circ f_{r-1}\circ f_r)$ equals $8T26$. The assertion of the lemma follows since the right-most indecomposable of $h$ is of degree $4$ or a prime.
			\end{proof}
			
		To conclude the proof of Theorem \ref{thm:main}, it remains to deal with right-tails $f_{r-3}\circ\dots\circ f_r$ of $f$ of degree sequence $(p,2,2,2)$ for $p$ prime, or $(4,2,2,2)$. That none of these right-tails can occur is the content of the following three lemmas.
			\begin{lemma}
			\label{lem:p222}
		   Right-tails of $f$ of degree sequence $(p,2,2,2)$ with an odd prime $p$ cannot occur.
		\end{lemma}
		
			\begin{lemma}
			\label{lem:4222}
			Right-tails of $f$ of degree sequence $(4,2,2,2)$ cannot occur.
		\end{lemma}
		
		\begin{lemma}
			\label{lem:2222}
			Right-tails of $f$ of degree sequence $(2,2,2,2)$ cannot occur.
			\end{lemma}
			
		The first of these three lemmas requires the most effort. We first deduce the other two.
		
		\begin{proof}[Proof of Lemmas \ref{lem:4222} and \ref{lem:2222} assuming Lemma \ref{lem:p222}]
			Another Magma search among subgroups of $C_2\wr C_2\wr C_2\wr S_4$ and of $C_2\wr C_2\wr C_2\wr C_2$ respectively, with the assumption $|\pi(K)|=2$, yields no candidate groups. But by Theorem \ref{thm:strong_diag}, $|\pi(K)|=4$ is impossible altogether for the tail degree sequence $(4,2,2,2)$, and for the sequence $(2,2,2,2)$ it implies the existence of an alternative decomposition of $f$ of tail degree sequence $(p,2,2,2)$ with $p$ odd. Since that one is assumed to be impossible for $p>1$, the only case left is $p=1$, i.e., $f$ of degree $16$, with monodromy group embedding into $C_2\wr C_2\wr C_2\wr C_2$.
			It now amounts to a straightforward check in Magma's database of transitive groups of degree $16$ to verify that no such $f$ exists.\footnote{In fact, this check does not require a search for genus-$0$ systems; the group-theoretical translation of minimal reducibility on the level of point and block stabilizers is enough to exclude any potential example.}
			\end{proof}
			
			The final piece required for the completion of the proof of Theorem \ref{thm:main} is therefore the proof of Lemma \ref{lem:p222}.
			\begin{proof}[Proof of Lemma \ref{lem:p222}]
				We may and will assume that $f_{r-3}\circ f_{r-2}$ admits no Ritt move, since otherwise we are reduced to the already excluded tail degree sequence $(p,2,2)$.
			Use the notation of Section \ref{sec:strong} with $s=4$, and note that the $2$-Sylow subgroup $P$ of $G_0\le \textrm{Mon}(f_{r-3}\circ f_{r-2}\circ f_{r-1}\circ f_r) \cap C_2\wr C_2\wr C_2\wr C_p$ is its intersection with the kernel $\ker(\textrm{Mon}(f_{r-3}\circ f_{r-2}\circ f_{r-1}\circ f_r)\to \textrm{Mon}(f_{r-3}))\le (C_2\wr C_2\wr C_2)^p$, which in particular is normal in $G_0$. It follows from Lemma \ref{lem:absorption} that $P/M$ (where $M=\pi(N_0\cap K)$) is of nilpotency class $\le 2=cl(C_2\wr C_2)$ and of exponent dividing $4$. 
			In other words, 
			\begin{equation}
				\label{eq:nilp}
				[[P,P],P]\le M
				\end{equation} and
				\begin{equation}
					 \label{eq:exp}
					 y^4\in M \text{ for all } y\in P.
					 \end{equation}
			
			Moreover, the map $\textrm{Mon}(f_{r-3}\circ\dots\circ f_{r})\to \textrm{Mon}(f_{r-3}\circ f_{r-2})$ maps $G_0$ onto $C_2\wr C_p$ (see \cite[Corollary 3.2(2)]{BKN26}) and thus maps $P$ onto $C_2^p$.
			Set $H:=\textrm{Mon}(f_{r-2}\circ f_{r-1}\circ f_r)\le C_2\wr C_2\wr C_2$, and let $H'\trianglelefteq H$ be the image of $P$ under projection to one block of size $8$ (in particular, $H'$ contains an $8$-cycle). We next note that $H'$
			satisfies the following Condition (Z).
			
			\begin{itemize}
				\item[(Z)] In the second center $Z_2(H')$ (i.e., the preimage of the center of $H'/Z(H')$ in $H'$), every element that interchanges the two blocks of size $4$ has order $8$.
				\end{itemize}

				This property in fact holds more generally, with the numbers $4$ and $8$ replaced by $2^{n-1}$ and $2^{n}$, $n\ge 3$,  for any subgroup  $H'\le S_{2^{n-1}}\wr C_2$ containing a cyclic transitive subgroup, but since we will soon require the precise shape of $H'$ anyway, it is most conveniently verified for our purposes by a direct check among imprimtive transitive subgroups $H'\le S_4\wr C_2$.
				
				We next show that Condition (Z) enforces $\textrm{rk}(K)\ge 3$ for all $p\ge 3$.
				Set $Y=Z_2(H')\cap H'_0$, where $H'_0\le H'$ denotes the subgroup of index $2$ preserving the two size-$4$ blocks. In particular, $P\cap Y^p$ maps  trivially to $C_2^p$. But since we have already seen that $P$ maps onto $C_2^p$, $P/(P\cap Y^p)$ has an elementary-abelian quotient of rank $p\ge 3$.  On the other hand, we know already from the proof of Lemma \ref{lem:red_to_4step} that $H'$ needs to be a normal subgroup of the group $8T17$ containing an $8$-cycle, and the only candidate groups $8T17$ and $8T7$ for such $H'$ do not have an elementary-abelian quotient of rank $\ge 3$.
				The surjection $\pi_1:P\to H'$ of $P$ to the first block therefore cannot induce an injection $P/(P\cap Y^p)\to \pi_1(P) Y/Y$. In particular, there must exist an element $x=(x_1,\dots, x_p)\in P$ whose first component lies in $Y$ whereas at least one other component does not lie in $Y$. Using $x$, we will construct an element $w\in M\le \pi(K)$ which is nontrivial on at least one, but not on all of the size-$8$ blocks. Let $S$ denote the set of components $i$ such that $x_i\notin Y$, and distinguish two cases. Firstly, if some $i\in S$ has $x_i\notin Z_2(H')$, then there are $y_i,z_i\in H'$ with $[[x_i,y_i],z_i]\ne 1$. Lift those to $y,z\in P$ and put $w=[[x,y],z]\in [[P,P],P]\le M$ (see \eqref{eq:nilp}). Then $w\ne 1$, but $w_j=1$ for all $j\notin S$ due to $x_j\in Z_2(H')$. Secondly, assume every $i\in S$ has $x_i\in Z_2(H')\setminus Y$. By definition of $Y$, such $x_i$ interchanges the size-$4$ blocks, and thus by Condition (Z) it has order $8$. But for $j\notin S$, $x_j\in Y\le H_0'$ has order dividing $4$. Thus $w=x^4\in M$ (see \eqref{eq:exp}) is as required.
				
				Now view $\pi(K)\le \mathbb{F}_2^{4p}$ as a module under the order-$p$ element $\tau=(\pi(\sigma^e)^8)$ (where $\langle\sigma\rangle$ is a cyclic transitive subgroup in $\textrm{Mon}(f)$, whence $\langle\pi(\sigma)\rangle$ is a cyclic transitive subgroup in $G_0$). 
				Since $\pi(K)$ is a semisimple $\langle\tau\rangle$-module, and also contains the diagonal submodule $\langle \pi(\sigma^e)^{4p}\rangle$, $\dim(\pi(K))\le 2$ implies that it is a sum of trivial modules, i.e., is centralized by $\tau$. But since $\tau$ permutes the size-$8$ blocks cyclically, the element $w\in \pi(K)$ constructed above is not fixed by $\tau$, showing $\dim(\pi(K))\ge 3$, contradicting Theorem \ref{thm:weak_diag}. 
				This concludes the proof of Lemma \ref{lem:p222}, and thus of (the forward direction of) Theorem \ref{thm:main}. For the converse, see the following section.
				\end{proof}
	
	\section{Explicit computations}
	Here we compute explicit equations for the minimally reducible pairs of polynomials $f$ and proper Siegel functions $g$ (over $\mathbb{Q}$) from Theorem \ref{thm:main}. We extract the key observation first, before presenting the concrete families:
	\begin{proposition}
		\label{thm:allpairs}
		Any minimally reducible pair $(f,g)$, where $f\in \mathbb{Q}[X]$ is a polynomial with solvable monodromy group and  $g\in \mathbb{Q}(X)$ is a proper Siegel function over $\mathbb{Q}$, belongs to one of the families computed in Theorems \ref{thm:8t44}, \ref{thm:8t28} and \ref{thm:8t47}.
		\end{proposition}
		\begin{proof}
			By Theorem \ref{thm:main}, all candidates have $\textrm{Mon}(f)\le S_8$, and a simple computer search (implementing the minimal reducibility condition on the level of point stabilizers) identifies all candidate monodromy groups together with the possible genus zero tuples (including in particular an $8$-cycle which acts as a product of two disjoint $8$-cycles in an additional degree-$16$ action). The output for the possible pairs $(G,A)$ of geometric and arithmetic monodromy groups yields only the following candidates:
			\begin{itemize}
				\item[1)] A class $4$-tuple of element orders $(2,2,2,8)$ in $G=A=C_2\wr S_4$, together with two class triples in the same group, arising from the $4$-tuple via coalescing. 
				See Theorem \ref{thm:8t44}.
				\item[2)] A class triple of element orders $(2,4,8)$ with $A=C_2\wr C_2\wr C_2$ and $[A:G]=2$. See Theorem \ref{thm:8t28}.
				\item[3)] A class $4$-tuple of orders $(2,2,2,8)$ in $G=A=S_4\wr C_2$, together with two class triples in the same group, arising from the $4$-tuple via coalescing. See Theorem \ref{thm:8t47}.
				\item[4)] A further class triple with $G$ of index $2$ in $S_4\wr C_2$. This is also a degeneration of the family in 3); however it cannot be defined over $\mathbb{Q}$; this can be verified via explicit computation, but has also been shown theoretically in \cite{Mue}; indeed, the candidate degree-$16$ action for the Siegel function would be primitive, implying an indecomposable Siegel function over $\mathbb{Q}$ with this geometric monodromy group, whose non-existence is however shown in \cite{Mue} (Proof of Theorem 5.2).
				\end{itemize}
				All that remains to do now is to explicitly parameterize the polynomials in  families 1), 2) and 3) up to equivalence (which is standard), and verify which parameters in these families indeed yield a Siegel function defined over $\mathbb{Q}$ as partner in the reducible pair. 
			\end{proof}
	\label{sec:ex}
	\subsection{Examples with monodromy group $C_2\wr S_4$}
	\begin{theorem}
		\label{thm:8t44}
	Let $s=ca^2-2cb^2-2$ where $c\in \mathbb{Q}\setminus\{0\}$, $a,b\in \mathbb{Q}$, $(a,b)\ne (0,0)$. Define the polynomial $f\in \mathbb{Q}[X]$ and the Siegel function $g\in \mathbb{Q}(X)$ as follows:
$$f_1 = X(X-4)(X+s)^2, \ f:=f_1(cX^2).$$
$$g_1 = -Y^2, \ g_2 = \frac{2(4(s+4)-Y^2)\cdot ((s+1)Y^2 + 4(s+2)Y + 4(s+4))}{(Y^2-8)^2},$$
$$g_3 =2(aY^2 + 4bY + 2a)/(bY^2+2aY+2b), \ g:=g_1\circ g_2\circ g_3.$$

Then $(f,g)$ is a minimally reducible pair. Conversely, every minimally reducible pair $(f,g)$ with $\textrm{Mon}(f)\cong C_2\wr S_4$ is of the above shape, up to simultaneous left- and right-composition of $f,g$ with linears.
\end{theorem}

\begin{proof}
	There are several ways to see that the given polynomials $f$ have monodromy group $C_2\wr S_4$ for all admissible values of $s$, the most convenient one being direct derivation from the group-theoretical data observed in Case 1) in the proof of Proposition \ref{thm:allpairs}, yielding a $4$-branch point family of polynomials $f$ with monodromy group $C_2\wr S_4$  in which exactly one simple preimage of a finite branch point of the $S_4$-quartic left-factor ramifies under the quadratic right-factor. From this, it is easy to see that any $4$-branch point member $f\in \mathbb{Q}[X]$ of this family must have the given shape, so far with the restriction $s\notin \{0,-4\}$ (but  yet without the restriction $s=ca^2-2cb^2-2$). The special values $s=0$ and $s=-4$ correspond to $3$-branch point polynomials with the same monodromy group.
	
	Note next that $g$ is indeed a Siegel function over the rationals: explicit computation shows that $g$ has denominator $(Y^2-2)^8$, i.e., it has exactly the two (real and algebraically conjugate) poles $\pm\sqrt{2}$, both of multiplicity $8$.
	One verifies directly that (the numerator of) $f(X)-g(Y)$ is reducible (indeed, it factors as a product of two polynomials of degree $4$ in $X$). 
	
	To show minimal reducibility, it is useful to note that $f = f_1(cX^2)$ is the unique decomposition of $f$ into indecomposable polynomials (up to the trivial manipulation $\tilde{f_1} = f_1\circ \mu$, $\tilde{f_2}=\mu^{-1}\circ cX^2$ for a linear polynomial $\mu$). Indeed, $f_1$ is indecomposable with primitive monodromy group $S_4$, and by Lemma \ref{lem:engstr}, polynomials of prime power degree always have a unique decomposition into indecomposables.
	As a matter of fact, the decomposition $g = g_1\circ g_2\circ g_3$ is also the unique decomposition of the Siegel function $g$ into indecomposable rational functions. Since this is however not quite as easy to see without digging into the group-theoretical details, we will refrain from using it as black box. 
	
	Instead we use the equality of splitting fields for minimally reducible pairs (Proposition \ref{prop:right_unique}a)). If $(f,g)$ were not minimal, then $g$ would either have to have a proper left-factor $\varphi$ forming a minimally reducible pair with $f_1$ (in particular, with the same splitting field of $f_1(X)-t$), or a proper left-factor $\varphi$ forming a minimally reducible pair with $f$. In the former case, since there are no nontrivial minimally reducible pairs with monodromy group $S_4$, we would have $\varphi=f_1$, and thus in particular the fiber product of $f_1$ and $g_1$ would have to be Siegel over $\mathbb{Q}$, since $g$ is. But that fiber product is birational to $X(X-4) = -Y^2$, and hence to $Y^2+U^2-4 = 0$, which has only finitely many rational points of bounded denominator, a contradiction. In the latter case, $\varphi$ would still have to have ramification index $8$ at infinity, while $\deg(\varphi)<16$, forcing $\varphi$ to be a polynomial. There are no minimally reducible pairs of two polynomials of degree $8$ (or indeed of any solvable polynomials of degree $>4$ by \cite[Theorem 6.1]{BKN26}), ruling out this possibility as well.

	To see that there are no other minimally reducible pairs with this monodromy group, it suffices to combine the computations from the proof of Proposition \ref{thm:allpairs} with a justification of the extra restriction $s=ca^2-2cb^2-2$. 
	This condition is obtained from parameterization of the genus-$0$ fixed field of the relevant index-$16$ subgroup, and corresponds exactly to the requirement of this fixed field being a rational function field, which is of course necessary to make $g$ Siegel over $\mathbb{Q}$. To make this process explicit, note that, as visible above, the octic left-factor $g_1\circ g_2$ of $g$ is always defined over $\mathbb{Q}$, independently of the value of $s$. Letting $v$ denote a root of $g_1(g_2(X))-t$, the pullback $f(X)-g_1(g_2(v))$ is still irreducible over $\mathbb{Q}(s,c)(v)$ (with $s,c$ taken as independent transcendentals), but reducible over $\mathbb{Q}(s,c)(v, \sqrt{(s+2)(v^2-8)/c})$, so that this quadratic extension must necessarily  correspond to the index-$16$ subgroup yielding the minimally reducible pair. For the Siegel function to be defined over $\mathbb{Q}$, we thus need $(s+2)(v^2-8)/c = w^2$ to have rational solutions $v,w$, and this condition is easily transformed into the condition $s=ca^2-2cb^2-2$ from the theorem statement.
	%
	\end{proof}

\begin{remark}
For completeness, we demonstrate how to obtain infinitely many explicit integral reducible values of bounded denominator for $f$ as in Theorem \ref{thm:8t44}. The process is analogous for the other families with a proper Siegel partner. For example, if $a,b,c\in \mathbb{Z}$, then $g$ is of the form $g(Y)=\frac{\gamma(Y)}{(Y^2-2)^8}$ for $\gamma(Y)\in \mathbb{Z}[Y]$. If $y\in \mathbb{Q}$ is chosen as ($\pm$) any partial convergent of the continued fraction expansion of $\sqrt{2}$, then $g(y)\in \mathbb{Z}$. To give a fully explicit example, choose $(a,b,c)=(1,0,1)$, so that $f=X^2(X^2-4)(X^2-1)^2$ and $g=-64Y^2(Y+1)^2(Y+2)^2(Y^2-2Y-2)^2(Y^2+2Y-2)^2/(Y^2-2)^8$. At the partial convergent $y=\frac{99}{70}$ of $\sqrt{2}$ one has $g(y)=-430164221433906960^2$ and $$f(X)-g(y) = (X^4 - 27720X^3 + 384199197X^2 - 18180690220838X + 430164221433906960)$$ $$\cdot (X^4 + 27720X^3 + 384199197X^2 + 18180690220838X + 430164221433906960).$$
 \end{remark}
 
 \subsection{Examples with monodromy group $C_2\wr C_2\wr C_2$}
 \begin{theorem}
 	\label{thm:8t28}
 	Let $a\in \mathbb{Q}$ be of the form $a=2u^2-v^2$ for $(u,v)\in \mathbb{Q}^2\setminus\{(0,0)\}$.
 	Let $$f_1=X^2, \ f_2 = X(X+4), f_3 = aX^2, \ f:=f_1\circ f_2\circ f_3 \in \mathbb{Q}[X].$$
 	and let
 %
 	$$g_1=-Y^2, \ g_2=2Y^2, \  
 	g_3 = \frac{2(3Y-8)(7Y-20)}{Y^2-8},\ 
 	g_4 = \frac{4(uY^2-2vY + 2u)}{-vY^2+4uY-2v},
 	\ g:=g_1\circ g_2\circ g_3\circ g_4.$$
 	Then $g\in \mathbb{Q}(X)$ is a Siegel function over $\mathbb{Q}$, and $(f,g)$ is a minimally reducible pair. Conversely, every minimally reducible pair $(f,g)$ with $\textrm{Mon}_{\mathbb{Q}}(f)\cong C_2\wr C_2\wr C_2$ is of the above shape, up to simultaneous left- and right-composition of $f,g$ with linears.
 	\end{theorem}

 	As in the previous subsection, the functions are computed using the explicit action of the monodromy group, and the minimal reducibility claim is again verifiable from the group-theoretical properties of this monodromy group, whence we refrain from going into the details again. Instead, it is worth pointing out some interesting properties of the pair in Theorem \ref{thm:8t28}. Firstly, while $f$ (and also $g$) has arithmetic monodromy group isomorphic to $C_2\wr C_2\wr C_2$, the geometric monodromy group is of index $2$ in this group (indeed, the quotient $\textrm{Mon}_{\mathbb{Q}}(g_1\circ g_2) = \textrm{Gal}(4Y^4+t / \mathbb{Q}(t))$ of $\textrm{Mon}_{\mathbb{Q}}(g)$ gives rise to a constant extension $\mathbb{Q}(i)/\mathbb{Q}$ inside the Galois closure), and inside this group, the group-theoretical equivalent of minimal reducibility is seen to no longer hold! In other words, $(f,g)$ is minimally reducible over $\mathbb{Q}$, but not over $\mathbb{Q}(i)$, which is quite obvious anyway from the fact that $g_1(X) = f_1(iX)$.
 	
 	 The group-theoretical interpretation of covers $\mathbb{P}^1\to \mathbb{P}^1$  in terms of branch cycle descriptions also yields an interesting connection between the pair in Theorem \ref{thm:8t28} and the family in Theorem \ref{thm:8t44}: Among the product-$1$  quadruples $(x_1,x_2,x_3,x_4)$ of cycle structure $([8], [2^3.1^2], [2^2.1^4], [2^2.1^4])$ generating the $C_2\wr S_4$ (corresponding to ramification types of a function $f$ as in Theorem \ref{thm:8t44}), there are some for which $(x_1, x_2x_3, x_4)$ is a triple of ramification type corresponding to the polynomial $f$ in Theorem \ref{thm:8t28}. This ``merging" of two inertia group generators corresponds topologically to the collision (via suitable choice of parameters) of two branch points in the family of Theorem \ref{thm:8t44}, so that,  up to linear equivalence (a priori, {\it over $\mathbb{C}$}), one should expect to regain the polynomial of Theorem \ref{thm:8t28} from this family via suitable specialization of parameters $a,b$. In fact, this would occur for the limit $s:=ca^2-2cb^2-2\to \infty$ in Theorem \ref{thm:8t44}, which requires some clarification and interpretation via projective closure and parameter change; namely,
 	 for $f_1:=X(X-4)(X+s)^2$, one has $\lim_{s\to\infty} \dfrac{2^8}{s^4} (f_1(cX^2))\circ (\sqrt{s}/2 X) = c^2X^4(cX^2+4)^2$, which is indeed the polynomial from Theorem \ref{thm:8t28} (with the parameter $a$ renamed as $c$).
 	 
 The distinct properties of the monodromy groups and thereby also of the involved polynomials observed throughout the paper justify considering the family in Theorem \ref{thm:8t28} as a separate type rather than simply a subfamily of the family in Theorem \ref{thm:8t44}.
 	 
 	 Another noteworthy property of the polynomial $f$ in Theorem \ref{thm:8t28}, observed in Case c) of Corollary \ref{cor:hilbert_siegel}, is that, due to its factoring through the biquadratic $f_1\circ f_2 = X^2(X+4)^2$ (which is equivalent up to composition with linears to a fourth degree Dickson polynomial, and hence itself part of a minimally reducible pair),  describing the exceptional reducible values of $f$ (up to finitely many) requires not just one, but {\it two} different Siegel functions; concretely, the function yielding infinitely many exceptional reducible values for the compositional left-factor $f_1\circ f_2$ of $f$ would be the polynomial $-4Y(Y+2)^2(Y+4) = h_1\circ h_2$ with 
 	 $h_1(Y) = -4Y^2+16$ and $h_2(Y)=Y^2+4Y+2$.
 	 
 	 Finally, note that the parameter value $a=-2$ in Theorem \ref{thm:8t28} allows for a particularly nice alternative parameterization up to composition with linear polynomials: In this case, one has $$\frac{1}{16} f - 1 = (X^2-1)^{\circ 3},$$
 	 the third iterate of $X^2-1$. In this form, the existence of the partner Siegel function forming the reducible pair and such leading to infinitely many reducible integral specialization values has been noticed in the context of arithmetic dynamics  (see \cite[Lemma 3.9]{Goksel} or \cite[Proposition 3.6]{IJOSS}). This is the special case $m=1$ of Theorem \ref{cor:newlyred}b).

 	 \subsection{Examples with monodromy group $S_4\wr C_2$}
 	 
 	 \begin{theorem}
 	 	\label{thm:8t47}
 	 	Let $s=a^2-2b^2-4$ with $a,b\in \mathbb{Q}$, $(a,b)\ne (0,0)$
 	 	and $$f(X) = (X^4 + 2sX^2 + (8s+ 32)X +s^2 -4s-24)^2.$$
 	 	Then there exists an indecomposable Siegel function $g(Y)\in \mathbb{Q}(Y)$ of degree $16$ such that $(f,g)$ is a minimally reducible pair. Conversely, $f$ belonging to a minimally reducible pair $(f,g)$ with $\textrm{Mon}(f)\cong S_4\wr C_2$ is of the above shape, up to  left- and right-composition  with linears.
 	 	\end{theorem}
 	 \begin{proof}
 	 	This is essentially due to M\"uller \cite{Mue}, see the proof of Theorem 5.2, which constructs the above family of polynomials $f$ (up to some trivial parameter changes) and notes that it shares a common Galois closure (with Galois group $S_4\wr C_2$) with an indecomposable Siegel function over $\mathbb{Q}$. As for the computation, note that generic values of $s$ correspond to a four-branch point family, with all finite branch points ramification index $2$, and with two of the three finite branch points of the ($S_4$) right-factor mapping to the same point under the left-factor $X^2$. From this, the computation of $f$ is standard. 
 	 	To complete the proof it only remains to verify that the (index-$16$) subgroup fixing a root of this Siegel function acts intransitively in the imprimitive degree-$8$ action of $S_4\wr C_2$, which is easily done and was in fact previously noted in \cite{BKN26}. Note that the restriction on $s$ once again arises due to requirement of the genus zero fixed field of the index-$16$ subgroup to be a rational function field. Note also that two particular admissible values of $s$, namely $s=-3$ and $s=-12$ (corresponding to $(a,b)=(1,0)$ and $(a,b)=(0,2)$) yield three-branch point polynomials with the same monodromy group $S_4\wr C_2$.
 	 	\end{proof}
 	 	\begin{remark}
 	 		The (indecomposable) Siegel function $g$ from Theorem \ref{thm:8t47} is somewhat complicated to write out; for purposes of verification, we give it in a form making use of certain symmetries between the coefficients.
 	 	Let 
	$$\gamma_1(Y) = \sum_{k=0}^3 c_k(a,b)Y^k  - c_k(-a,-b)(\frac{2}{Y})^k, \  \gamma_2(Y) = \sum_{k=0}^2 d_k(a,b) Y^k + d_k(-a,-b)(\frac{2}{Y})^k,$$ 
 	 	where $$c_0=4a^3 + 16a^2b + 12ab^2 - 12a + 8b^3 - 16b,$$
 	 	$$c_1 = 4a^3 + 12a^2b + 6a^2 + 22ab^2 + 30ab - 4a + 6b^3 + 18b^2 - 6b - 6,$$
 	 	$$c_2 = 4a^2b + 6a^2 + 6ab^2 + 12ab + 6a + 4b^3 + 12b^2 + 8b,$$
 	 	$$c_3 = a^2 + ab^2 + 3ab + 2a + b^3 + b^2 + 3b + 1,$$
 	 	$$d_0  = 4a^2 + 12ab + 10b^2 + 2,$$
 	 	$$d_1 = 4a^2 + 12ab + 12a + 8b^2 + 16b,$$
 	 	$$d_2 = a^2 + 2ab + 4a + b^2 + 6b - 1.$$
 	 	Then $g$ is defined as $$g(Y) = \frac{-2^6 \gamma_1^2 \gamma_2 }{(Y-\frac{2}{Y})^8}.$$
 	 		\end{remark}
 	 		
 	 		A final remark on the (non-)uniqueness of the partners in the pairs in Theorems \ref{thm:8t44}, \ref{thm:8t28}, \ref{thm:8t47}.
 	 		\begin{remark}
 	 			\label{rem:2to1}
 	 			It can be verified purely group-theoretically with Magma (searching for suitable genus-$0$ tuples in permutation actions corresponding to a minimally reducible pair) that each of the polynomials $f$ in Theorems \ref{thm:8t44}, \ref{thm:8t28}, \ref{thm:8t47} has a unique (up to right-composition with linears) Siegel partner in a minimally reducible pair. This is indeed important for the counting of exceptional partners in Corollary \ref{cor:hilbert_siegel}. 
 	 			Conversely, the same kind of Magma search yields that, for the functions in Theorem \ref{thm:8t47} each $g$ also has a unique partner $f$, whereas in Theorems \ref{thm:8t44} and \ref{thm:8t28}, the mapping from $f$ to $g$ is exactly $2$-to-$1$! Concretely, in Theorem \ref{thm:8t44}, $f_1(cX^2)$ and $f_1(4-cX^2)$ share the same $g$ (note that $f_1(4-cX^2)$ is admissible, since $f_1(4-X)=X(X-4)(X-(s+4))^2$ is an admissible form for the left-factor; indeed, if $s=ca^2-2cb^2-2$, then $-(s+4)=-ca^2+2cb^2-2 = c(a+2b)^2-2c(a+b)^2-2$, which is an admissible value for the same parameter $c$). Similarly, in Theorem \ref{thm:8t28}, swapping the parameter $a$ for $-a$ yields the same $g$.
 	 			\end{remark}

\end{document}